\documentclass[12pt]{amsart}
\usepackage{amsmath, amsfonts, amssymb, amsthm, eufrak}
\usepackage[all]{xy}
\newtheorem{sub}{}[section]
\newtheorem{subsub}{}[sub]
\let\a\alpha
\let\b\beta

\let\z\zeta

\let\l\lambda
\let\m\mu

\let\s\varsigma

\let\f\varphi
\let\o\omega

\let\L\Lambda

\let\Om\Omega

\def\mM#1{{\hbox{$#1$}}}
\def\lra{\longrightarrow}

\def\C{\mathbb C}
\def\CC{\mathcal C}

\def\E{\mathcal E}
\def\EE{\text{E}}
\def\F{\mathcal F}
\def\G{\mathcal G}
\def\H{\text{H}}
\def\I{\mathcal I}

\def\L{\mathcal L}
\def\N{\mathcal N}
\def\M{\text{M}_{\P^2}}
\def\MM{\text{\emph{M}}_{\P^2}}
\def\P{\mathbb{P}}
\def\W{\mathbb{W}}
\def\R{\text{R}}

\def\T{\mathcal T}

\def\O{\mathcal O}

\def\q{\mathfrak{q}}  

\def\Ker{{\mathcal Ker}}
\def\Coker{\mathcal Coker}
\def\Im{{\mathcal Im}}
\def\Hom{\text{Hom}}
\def\Aut{\text{Aut}}
\def\Ext{\text{Ext}}
\def\GL{\text{GL}}

\def\sigg{\mathop{\hbox{$\displaystyle\sum$}}\limits}
\def\nsp{\lbrace 0\rbrace}
\def\Ssect#1#2{\pagebreak[3]\begin{sub}\label{#2}{\sc\small\small
#1}\rm\medskip}
\def\ot{\otimes}

\def\tensor{\otimes}
\def\isom{\simeq}

\newcommand{\tilda}{\widetilde}

\def\bdm{\begin{displaymath}}
\def\edm{\end{displaymath}}

\def\ba{\begin{array}}
\def\ea{\end{array}}
\def\sepsec{\vskip 1.4cm}
\def\sepsub{\vskip 1cm}
\def\SEP{\vskip 0.8cm}
\def\sEP{\vskip 0.6cm}
\def\carre{$\Box$}
\def\flinc{\ar@{^{(}->}}
\def\fleq{\ar@{=}}
\def\flon{\ar@{->>}}
\def\fmaps{\ar@{|-{>}}}

\def\Proof{\noindent\textsc{Proof:}\ }

\begin{document}

\title[moduli spaces of semistable sheaves on plane quartics]
{on the geometry of the
moduli spaces of semistable sheaves supported on plane quartics}

\author{Jean-Marc Dr\'ezet \and Mario Maican}

\address{
Jean-Marc Dr\'ezet \\
Institut de Math\'ematiques \\
173 Rue du Chevaleret \\
F-75013 Paris, France \\
E-mail: drezet@math.jussieu.fr}
\address{
Mario Maican \\
1 University Circle \\
Western Illinois University \\
Macomb, IL 61455, USA \\
E-mail: m-maican@wiu.edu}

\begin{abstract}
We decompose each moduli space of semistable sheaves on the complex projective
plane with support of dimension one and degree four into locally closed
subvarieties, each subvariety being the good or geometric quotient of a set of
morphisms of locally free sheaves modulo a reductive or a nonreductive group.
We find locally free resolutions of length one of all these sheaves and
describe them.
\end{abstract}

\maketitle

\tableofcontents


\section{Introduction}

Let $V$ be a three dimensional vector space over $\C$, and $\P^2=\P(V)$ the
projective plane of lines in $V$. Let $\M(r,\chi)$ denote the moduli space of
semistable sheaves $\F$ on $\P^2$ with Hilbert polynomial $P_{\F}(t)=rt+\chi$.
The positive integer $r$ is the multiplicity of $\F$ while $\chi$ is its Euler
characteristic. The generic stable sheaves in this moduli space are the line
bundles of Euler characteristic $\chi$ on smooth plane curves of degree $r$. The
map sending $\F$ to the twisted sheaf $\F(1)$ gives an isomorphism between
$\M(r,\chi)$ and $\M(r,r+\chi)$, so we can restrict our attention to the case $0
< \chi \le r$. It is known from \cite{lepotier} that the spaces $\M(r,\chi)$ are
projective, irreducible, locally factorial, of dimension $r^2 +1$, and smooth at
all points given by stable sheaves.

It is easy to see that $\M(2,1)$ is isomorphic to the space of conics in $\P^2$
while $\M(2,2)$ is the good quotient modulo the action by conjugation
of the group $\GL(2,\C) \times \GL(2,\C)$ on the space of $2 \times
2$-matrices with entries in $V^*$ and nonzero determinant.

The case of multiplicity three is also well-understood. J.~Le Potier showed in
\cite{lepotier} that $\M(3,2)$ and $\M(3,1)$ are both isomorphic to the
universal cubic in $\P(V) \times \P(S^3V^*)$.
It was first noticed in \cite{maican-diplom} that $\M(3,2)$
is the geometric quotient of the set of injective morphisms \ \mM{\f:
\O(-2) \oplus \O(-1)\to 2\O} \ for which $\f_{12}$ and $\f_{22}$ are linearly
independent regarded as elements of $V^*$, modulo the action by conjugation of
the nonreductive algebraic group \ $\Aut(\O(-2) \oplus \O(-1)) \times
\Aut(2\O)$. The corresponding result for M$_{\P^2}(3,1)$ was established in
\cite{freiermuth-diplom} and \cite{freiermuth-trautmann}.

According to \cite{lepotier}
$\M(3,3)$ contains an open dense subset which is
a good quotient of the space of injective morphisms
$3 \O(-1) \to 3\O$ modulo the action by conjugation of
$\GL(3,\C)\times \GL(3,\C)$. The complement of this set is
isomorphic to the space of cubics $\P(S^3 V^*)$.
We summarize in the following table the facts about $\M(r,\chi)$, $r=1,2,3$,
that are known:

\begin{center}
\noindent \\
\begin{tabular}{|l|p{8cm}|}
\hline
\multicolumn{2}{|l|}{$\M(1,1)$} \\
\hline
\begin{tabular}{l}
$h^0(\F(-1))=0$ \\
$h^1(\F)=0$\\
$h^0(\F \tensor \Om^1(1))=1$
\end{tabular}
&
$0 \lra \O(-1) \lra \O \lra \F \lra 0$ \\
\hline
\end{tabular}
\end{center}
\begin{center}
\begin{tabular}{|l|p{8cm}|}
 \hline
\multicolumn{2}{|l|}{$\M(2,1)$} \\
\hline
\begin{tabular}{l}
$h^0(\F(-1))=0$ \\
$h^1(\F)=0$\\
$h^0(\F \tensor \Om^1(1))=0$
\end{tabular}
&
$ 0 \lra \O(-2) \lra \O \lra \F \lra 0$ \\
\hline
\end{tabular}
\end{center}
\begin{center}
\begin{tabular}{|l|p{8cm}|}
\hline
\multicolumn{2}{|l|}{$\M(2,2)$} \\
\hline
\begin{tabular}{l}
$h^0(\F(-1))=0$ \\
$h^1(\F)=0$\\
$h^0(\F \tensor \Om^1(1))=2$
\end{tabular}
&
$0 \lra 2\O(-1) \lra 2\O \lra \F \lra 0$ \\
\hline
\end{tabular}
\end{center}
\begin{center}
\begin{tabular}{|l|p{8cm}|}
\hline
\multicolumn{2}{|l|}{$\M(3,1)$} \\
\hline
\begin{tabular}{l}
$h^0(\F(-1))=0$ \\
$h^1(\F)=0$\\
$h^0(\F \tensor \Om^1(1))=0$
\end{tabular}
&
$0 \lra 2\O(-2) \lra \O(-1) \oplus \O \lra \F \lra 0$ \\
\hline
\end{tabular}
\end{center}
\begin{center}
\begin{tabular}{|l|p{8cm}|}
\hline
\multicolumn{2}{|l|}{$\M(3,2)$} \\
\hline
\begin{tabular}{l}
$h^0(\F(-1))=0$ \\
$h^1(\F)=0$ \\
$h^0(\F \tensor \Om^1(1))=1$
\end{tabular}
&
$0 \lra \O(-2) \oplus \O(-1) \lra 2\O \lra \F \lra 0$ \\
\hline
\end{tabular}
\end{center}
\begin{center}
\begin{tabular}{|l|p{8cm}|}
\hline
\multicolumn{2}{|l|}{$\M(3,3)$} \\
\hline
\begin{tabular}{l}
$h^0(\F(-1))=0$ \\
$h^1(\F)=0$ \\
$h^0(\F \tensor \Om^1(1))=3$
\end{tabular}
&
$0 \lra 3\O(-1) \lra 3\O \lra \F \lra 0$ \\
\hline
\begin{tabular}{l}
$h^0(\F(-1))=1$ \\
$h^1(\F)=0$ \\
$h^0(\F \tensor \Om^1(1))=3$
\end{tabular}
&
$0 \lra \O(-2) \lra \O(1) \lra \F \lra 0$ \\
\hline
\end{tabular}

\end{center}

\medskip

For each moduli space (except $\M(3,3)$) the left column indicates the
cohomological conditions verified by the corresponding sheaves. These sheaves
are isomorphic to cokernels of morphisms of locally free sheaves described in
the right column. The moduli space is isomorphic to the good quotient, modulo
the action of the appropriate group, of a certain open subset of the set of
these morphisms. The moduli space $\M(3,3)$ is the disjoint union of a dense
open subset described in the first line, and of a locally closed subset
described in the second line.

\medskip

In this paper we will study the spaces $\M(4,\chi)$ for $1\leq\chi\leq 4$..
We will decompose each moduli space into locally closed subvarieties (which we
call {\em strata}) given by cohomological conditions and we will describe these
subvarieties as good or geometric quotients of spaces of morphisms. The
work of finding resolutions for sheaves $\F$ in $\M(4,\chi)$, apart from the
case $\chi=4$, $h^0(\F(-1))=1$, has already been carried out in \cite{maican}
and is summarized in the next table. Each stratum in $\M(4,\chi)$
described by the cohomological conditions from the left column of the table
below is isomorphic to the good quotient, modulo the action of the appropriate
group, of a certain open subset of the set of morphisms of locally free  sheaves
from the middle column. The right column gives the codimension of the stratum.

\begin{center}
\noindent \\
\begin{tabular}{|l|p{9cm}|c|}
\hline
\multicolumn{3}{|l|}{$\M(4,1)$} \\
\hline
\begin{tabular}{l}
$h^0(\F(-1))=0$ \\
$h^1(\F)=0$\\
$h^0(\F \tensor \Om^1(1))=0$
\end{tabular}
&
$0 \to 3\O(-2) \to 2\O(-1) \oplus \O \to \F \to 0$ & 0\\
\hline
\end{tabular}
\end{center}
\begin{center}
\begin{tabular}{|l|p{9cm}|c|}
\hline
\begin{tabular}{l}
$h^0(\F(-1))=0$ \\
$h^1(\F)=1$\\
$h^0(\F \tensor \Om^1(1))=1$
\end{tabular}
&
$0 \to \O(-3) \oplus \O(-1) \to 2\O \to \F \to 0$ & 2\\
\hline
\end{tabular}
\end{center}
\begin{center}
\begin{tabular}{|l|p{9cm}|c|}
\hline
\multicolumn{3}{|l|}{$\M(4,2)$} \\
\hline
\begin{tabular}{l}
$h^0(\F(-1))=0$ \\
$h^1(\F)=0$\\
$h^0(\F \tensor \Om^1(1))=0$
\end{tabular}
&
$0 \to 2\O(-2) \to 2\O \to \F \to 0$ & 0\\
\hline
\begin{tabular}{l}
$h^0(\F(-1))=0$ \\
$h^1(\F)=0$\\
$h^0(\F \tensor \Om^1(1))=1$
\end{tabular}
&
$0 \to 2\O(-2) \oplus \O(-1) \to \O(-1) \oplus 2\O \to \F \to 0$ & 1\\
\hline
\begin{tabular}{l}
$h^0(\F(-1))=1$ \\
$h^1(\F)=1$\\
$h^0(\F \tensor \Om^1(1))=3$
\end{tabular}
&
$0 \to \O(-3) \to \O(1) \to \F \to 0$ & 3\\
\hline
\end{tabular}
\end{center}
\begin{center}
\begin{tabular}{|l|p{9cm}|c|}
\hline
\multicolumn{3}{|l|}{$\M(4,3)$} \\
\hline
\begin{tabular}{l}
$h^0(\F(-1))=0$ \\
$h^1(\F)=0$\\
$h^0(\F \tensor \Om^1(1))=2$
\end{tabular}
&
$0 \to \O(-2) \oplus 2\O(-1) \to 3\O \to \F \to 0$ & 0\\
\hline
\begin{tabular}{l}
$h^0(\F(-1))=1$ \\
$h^1(\F)=0$\\
$h^0(\F \tensor \Om^1(1))=3$
\end{tabular}
&
$0 \to 2\O(-2) \to \O(-1) \oplus \O(1) \to \F \to 0$ & 2\\
\hline
\end{tabular}
\end{center}
\begin{center}
\begin{tabular}{|l|p{9cm}|c|}
\hline
\multicolumn{3}{|l|}{$\M(4,4)$} \\
\hline
\begin{tabular}{l}
$h^0(\F(-1))=0$ \\
$h^1(\F)=0$\\
$h^0(\F \tensor \Om^1(1))=4$
\end{tabular}
&
$0 \to 4\O(-1) \to 4\O \to \F \to 0$ & 0\\
\hline
\begin{tabular}{l}
$h^0(\F(-1))=1$ \\
$h^1(\F)=0$\\
$h^0(\F \tensor \Om^1(1))=4$
\end{tabular}
&
$0 \to \O(-2) \oplus \O(-1) \to \O \oplus \O(1) \to \F \to 0$ & 1\\
\hline
\end{tabular}
\end{center}

\bigskip

One of our difficulties will be to show that the natural maps from the spaces of
morphisms of sheaves to the subsets in the moduli spaces are good or geometric
quotient maps. At \ref{4.E} and \ref{4.E3} the difficulty is compounded by the
fact that the group is not reductive and we do not know a priori the existence
of a good quotient. A similar situation was considered by the first author in
\cite{drezet-1991}, cf. the proof of theorem D. Our method, exhibited in the
proof of \ref{3.B}, is reminiscent of \cite{drezet-1991}, in that we use an
already established quotient modulo a reductive group, but different, in that we
do not have a diagram as at 5.4 in \cite{drezet-1991}, but only local morphisms.
These local morphisms are constructed using the relative Beilinson spectral
sequence.

The second chapter of this paper contains some tools that are used in the next
chapters to study the moduli spaces. The chapters 3,4,5 are devoted to
$\M(4,1)$ and $\M(4,3)$, $\M(4,2)$, $\M(4,4)$ respectively. 
Now we summarize the descriptions and properties of the moduli spaces
$\M(4,\chi)$ that are contained in this paper :
\sepsub

\Ssect{The moduli spaces $\M(4,1)$ and $\M(4,3)$}{m4143}

These moduli spaces are isomorphic : $\M(4,1)\simeq\M(4,-1)$ by duality (cf.
\ref{dual}) and $\M(4,-1)\simeq\M(4,3)$ (the isomorphism sending to the point
corresponding to the sheaf $\F$ the point corresponding to $\F(1)$). So we will
treat together in chapter \ref{M4143} the moduli spaces $\M(4,1)$ and $\M(4,3)$.
These varieties are smooth because in this case semi-stability is equivalent to
stability. We show that, as in the preceeding table, $\M(4,3)$ is the disjoint
union of two strata : $\M(4,3)=X_0(4,3)\amalg X_1(4,3)$, where $X_0(4,3)$ is an
open subset and $X_1(4,3)$ a closed smooth subvariety of codimension 2. These
strata correspond to sheaves having the cohomological conditions given in the
first column of the preceeding table. For example the open stratum $X_0(4,3)$
contains the points representing the sheaves $\F$ such that \
$h^0(\F(-1))=h^1(\F)=0$, $h^0(\F\otimes\Omega^1(1))=2$.

The sheaves in the open stratum $X_0(4,3)$ are isomorphic to cokernels of
injective morphisms
$$\O(-2)\oplus 2\O(-1)\lra 3\O \ \ .$$
Let $W=\text{Hom}(\O(-2)\oplus 2\O(-1),3\O)$, on which acts the non reductive
algebraic group
$$G=\big(\Aut(\O(-2)\oplus 2\O(-1))\times\Aut(3\O)\big)/\C^*$$
in an obvious way. Following \cite{drezet-trautmann} we describe in \ref{MSM}
and \ref{open1} an open $G$--invariant subset $\mathbf{W}$ of $W$ such that
there exists a geometric quotient $\mathbf{W}/G$ which is a smooth projective
variety. Then $X_0(4,3)$ is canonically isomorphic to the open subset of
$\mathbf{W}/G$ corresponding to injective morphisms.

The sheaves in the closed stratum $X_1(4,3)$ are isomorphic to cokernels of
injective morphisms
$$2\O(-2)\lra\O(-1)\oplus\O(1) \ ,$$
and $X_1(4,3)$ is isomorphic to a geometric quotient of a suitable open subset
of \hfil\break$\Hom(2\O(-2), \O(-1)\oplus\O(1))$ by the non reductive group
$$\big(\Aut(2\O(-2))\times\Aut(\O(-1)\oplus\O(1))\big)/\C^* \ .$$

Similarly $\M(4,1)$ is the disjoint union of the open subset $X_0(4,1)$ and the
closed smooth subvariety $X_1(4,1)$. Of course the canonical isomorphism
$\M(4,1)\simeq \M(4,3)$ induces isomorphisms of the strata.

We can give a precise description of the sheaves that appear in the strata. We
define a closed subvariety $\widetilde{Y}$ of $X_0(4,3)$ corresponding to
sheaves $\F$ such that there is a non trivial extension
\[0\lra\O_\ell(-1)\lra\F\lra\O_X(1)\lra 0\]
where $\ell$ is a line and $X$ a cubic. The open subset
$X_0(4,3)\backslash\widetilde{Y}$ consists of kernels of surjective morphisms
$\O_C(2)\to\O_Z$, where $C$ is a quartic and $Z$ a length 3 finite subscheme of
$\P^2$ not contained in any line. The sheaves in $X_1(4,1)$ are the kernels of
the surjective morphisms $\O_C(1)\to\O_P$, $C$ beeing a quartic and $P$ a
closed point of $C$. The sheaves in the other strata can be described similarly
using the isomorphism $\M(4,1)\simeq \M(4,3)$.
\end{sub}

\sepsub

\Ssect{The moduli space $\M(4,2)$}{m42}

It is treated in chapter 4.
Here we have 3 strata : an open subset $X_0$, a locally closed smooth
subvariety $X_1$ of codimension 1, and a closed smooth subvariety
$X_2$ of codimension 3 which is in the closure of $X_1$. These strata are
defined by cohomological conditions on the corresponding sheaves, as indicated
in the preceeding table.

The smallest stratum $X_2$ is the set of sheaves of the form $\O_C(1)$, where
$C$ is a quartic. Hence $X_2$ is isomorphic to $\P(S^4V^*)$.

The stratum $X_0$ contains the cokernels of the injective semi-stable morphisms
$$2\O(-2)\lra 2\O \ .$$
Let $N(6,2,2)$ denote the moduli space of semi-stable morphisms $2\O(-2)\to 2\O$
(cf. \ref{KM}). The closed subset corresponding to non injective morphisms is
naturally identified with $\P^2\times\P^2$, hence $X_0$ is isomorphic to
$N(6,2,2)\backslash(\P^2\times\P^2)$.

The sheaves of the second stratum $X_1$ are the cokernels of injective morphisms
$$2\O(-2)\oplus\O(-1)\lra\O(-1)\oplus 2\O$$
corresponding to matrices
$\f= \begin{pmatrix}
X_1 & X_2 & 0 \\
\star & \star & Y_1 \\
\star & \star & Y_2
\end{pmatrix}$
where $X_1,X_2 \in V^*$ are linearly independent one-forms and the same for
$Y_1,Y_2 \in V^*$.
Let
$$G=\Aut(2\O(-2) \oplus \O(-1)) \times \Aut(\O(-1) \oplus 2\O) \ .$$
This non reductive group acts naturally on the variety $W_1$ of the preceeding
matrices and we prove that there is a geometric quotient $W_1/G$ which is
isomorphic to $X_1$. Using this description we conclude that the generic sheaf
in $X_1$ is of the form $\O_C(1)(P-Q)$, where $C$ is a smooth quartic and $P$,
$Q$ are distinct points of $C$.

Let $X=X_0\cup X_1$. We prove that the sheaves in $X$ are precisely the
cokernels of the injective morphisms
$$2\O(-2)\oplus\O(-1)\lra\O(-1)\oplus 2\O$$
which are $G$--semi-stable with respect to the
polarization $(\frac{1-\mu}{2},\mu,\mu,\frac{1-\mu}{2})$, where $\mu$ is a
rational number such that $\frac{1}{3}<\mu<\frac{1}{2}$ (cf.
\cite{drezet-trautmann}, \ref{MSM}). Let $W$ be the set of such injective
$G$--semi-stable morphisms. We prove that there is a good quotient $W//G$ which
is isomorphic to $X$. The existence of this quotient cannot be obtained from
the results of \cite{drezet-trautmann} and \cite{dr3b}.

The inclusion $X_0\subset N(6,2,2)$ can be extended to a morphism \
$\delta:X\to N(6,2,2)$. Let $\widetilde{\mathbf{N}}$ denote the blowing-up of
$N(6,2,2)$ along $\P^2\times\P^2$. We prove that $X$ is naturally isomorphic to
an open subset of $\widetilde{\mathbf{N}}$ and $\delta$ can be identified to
the restriction to $X$ of the natural projection $\widetilde{\mathbf{N}}\to
N(6,2,2)$. The complement of $X$ in $\widetilde{\mathbf{N}}$ is contained in
the inverse image of the diagonal of $\P^2\times\P^2\subset N(6,2,2)$, and for
every $P\in\P^2$ the inverse image of $(P,P)$ contains exactly one point of
$\widetilde{\mathbf{N}}\backslash X$. Hence $\widetilde{\mathbf{N}}\backslash X$
is isomorphic to $\P^2$, whereas $\M(4,2)\backslash X$ is isomorphic to
$\P(S^4V^*)$.
\end{sub}

\sepsub

\Ssect{The moduli space $\M(4,4)$}{m44}

It has been completely described by J. Le Potier in \cite{lepotier}.
Let $N(3,4,4)$ be the moduli space of semi-stable morphisms $4\O(-1)\to 4\O$
(cf. \ref{KM}). Then $\M(4,4)$ is isomorphic to the blowing-up of $N(3,4,4)$
along a subvariety isomorphic to $\P(V^*)$.

The open stratum $X_0$ of $\M(4,4)$ contains the sheaves which are cokernels of
injective morphisms $4\O(-1)\to 4\O$ (which are then semi-stable). Our
contribution to the study of $\M(4,4)$ is a description of the complement \
$X_1=\M(4,4)\backslash X_0$, i.e of the exceptional divisor in $\M(4,4)$. We
show that the sheaves in $X_1$ are precisely the cokernels of the injective
morphisms
$$\O(-2)\oplus\O(-1)\lra\O\oplus\O(1)$$
such that $\phi_{12}\not=0$. Let \ $W=\Hom(\O(-2)\oplus\O(-1),\O\oplus\O(1))$,
on which acts the non reductive group
$$G=\big(\Aut(\O(-2)\oplus\O(-1))\times\Aut(\O\oplus\O(1))\big) \ .$$
les $W_1$ be the set of $G$-semi-stable points of $W$ with respect to a
polarization $(\lambda_1,\lambda_2,\mu_1,\mu_2)$ such that
$\lambda_1=\mu_2<\frac{1}{4}$ (cf. \ref{MSM}). According to
\cite{drezet-trautmann} there exists a good quotient $W_1//G$ which is a
projective variety. We prove that $X_1$ and $W_1//G$ are canonically isomorphic.

\end{sub}

\sepsub

\Ssect{Clifford's theorem}{clif}

To show that there are no sheaves in $\M(4,\chi)$ other than those from the
table we prove at \ref{3.E}, \ref{4.B} and \ref{5.F} the following cohomology
estimates:

\noindent \\
\emph{For any sheaf $\F$ in $\MM(4,0)$ or in $\MM(4,1)$ we have $h^1(\F) \le
1$.}

\noindent \\
\emph{For any sheaf $\F$ in $\MM(4,2)$ we have $h^1(\F)=0$, unless $\F \isom
\O_C(1)$ for a quartic $C \subset \P^2$.}

\noindent \\
For generic sheaves in $\M(4,\chi)$ the above estimates already follow from
Clifford's theorem. Indeed, a generic sheaf in $\M(r,\chi)$ is a line bundle
supported on a smooth curve of degree $r$. Clifford's theorem, cf.
\cite{griffiths-harris} p. 251, states that, if $\L$ is a line bundle on a
compact Riemann surface $S$ corresponding to an effective divisor and such that
$h^1(\L)>0$, then we have the inequality
\bdm
h^0(\L) \le 1 + \frac{\text{deg}(\L)}{2},
\edm
with equality only if $\L = \O_S$ or $\L = \o_S$ or $S$ is hyperelliptic.

\noindent \\
If $\L$ has Hilbert polynomial $P(t)=rt+\chi$, the Riemann-Roch theorem and the
genus formula give
\bdm
\text{deg}(\L)= g(S)-1+\chi = \frac{(r-1)(r-2)}{2}-1+\chi = \frac{r(r-3)}{2} +
\chi
\edm
and the above inequality takes the form
\bdm
h^0(\L) \le 1 + \frac{\chi}{2} + \frac{r(r-3)}{4}.
\edm
Taking $r=4$ and $0 \le \chi < 4$ and noting that equality in Clifford's 
theorem is achieved for non-generic sheaves, we conclude that for a generic
sheaf $\F$ in $\M(4,\chi)$ we have the relation
\bdm
h^0(\F) < 2 + \frac{\chi}{2}.
\edm
This yields \ref{3.E}, \ref{4.B} and \ref{5.F} for generic sheaves. What we
prove in this paper is that, in fact, Clifford's theorem is true for all sheaves
in $\M(4,\chi)$, $0\leq\chi<4$. By inspecting the first table from above we see
that Clifford's theorem is also true for all sheaves in $\M(r,\chi)$, $r=1,2,3$,
$0\leq \chi<r$. There is thus enough evidence to suggest that the following
``generalized Clifford's theorem" be true:

\noindent \\
{\bf Conjecture:} {\em Let $\F$ be a semistable sheaf on $\P^2$ with Hilbert
polynomial $P(t)=rt+\chi$, where $r \ge 1$ and $0\leq\chi<r$ are integers.
If $h^1(\F)>0$, then we have the inequality
\bdm
h^0(\F) \ \leq \ 1 + \frac{\chi}{2} + \frac{r(r-3)}{4}
\edm
with equality only in the following two cases :

(i) $r=3$, $\chi=0$ and $\F=\O_C$ for some cubic $C$.

(ii) $r=4$, $\chi=2$ and $\F=\O_C(1)$ for some quartic $C$.}
\end{sub}

\sepsec

\section{Preliminaries}

\Ssect{Duality}{dual}

For a sheaf $\F$ on $\P^2$ with support of dimension one we will consider the
dual sheaf \hfil\break$\F^D={\mathcal Ext}^1(\F,\omega_{\P^2})$. According to
\cite{maican-duality} we have the following duality result:

\SEP

\begin{subsub}\label{2.A}{\bf Theorem : }
The map $\xymatrix{\F \fmaps[r]& \F^D}$ gives an isomorphism \hfil\break
$\MM(r,\chi) \stackrel{\simeq}{\lra}\MM(r,-\chi)$.
\end{subsub}

\SEP

For a semistable sheaf $\F$ on $\P^2$ with support of dimension one Serre 
duality gives the relations

\SEP

\begin{subsub}\label{2.B}{\bf Proposition : }
We have \ $h^i(\F \tensor \Om^j(j))=h^{1-i}(\F^D \tensor
\Om^{2-j}(3-j))$, in particular \ $h^i(\F)=h^{1-i}(\F^D)$.
\end{subsub}

\SEP

\begin{subsub}\label{2.C}{\bf Proposition : }
Let $\F$ be a semistable sheaf on $\P^2$ with Hilbert polynomial\hfil\break
$P_\F(t)=rt+\chi$. Then we have \ $h^0(\F(i))=0$ \ for \ ${\displaystyle i <
\frac{3-r}{2}- \frac{\chi}{r}}$ , and \ $h^1(\F(i))=0$ \ for \ ${\displaystyle i
> \frac{r-3}{2}-\frac{\chi}{r}}$.
\end{subsub}

\Proof
Let $D$ be the schematic support of $\F$ (defined by the associated Fitting
ideal). It is a curve of degree $r$ and $\F$ can be viewed as a sheaf on $D$.
Assume that $h^0(\F(i))>0$. Then there is a nonzero morphism $\O\to\F(i)$,
inducing a nonzero morphism $\sigma:\O_D\to\F$. Let $\I=\ker(\sigma)$. We
claim that $\sigma$ factors through an injective morphism $\O_C \to \F(i)$,
for a curve $C\subset\P^2$ contained in $D$. Let $f\in H^0(\O(r))$ be an
equation of $C$. According to 6.7 in \cite{maican}, there is a polynomial
$g$ dividing $f$ such that $\I$ is contained in the ideal sheaf $\G$ in $\O_D$
defined by $g$ and that $\G/\I$ is supported on finitely many points. Since
$\F(i)$ has no zero dimensional torsion, $\G/\I$ is mapped to zero in $\F(i)$,
hence $\G=\I$ and we may take for $C$ the support of $\G$.

Let $d=\deg(C)\leq r$. From the semistability of $\F(i)$ we get
\bdm
\frac{3-d}{2} = p(\O_C) \ \le \ p(\F) = \frac{\chi}{r} +i, \quad \text{hence}
\quad
\frac{3-r}{2} - \frac{\chi}{r} \le i.
\edm
This proves the first part of proposition \ref{2.C}. The second part follows
from the first and \ref{2.B}, \ref{2.A}.
\carre

\end{sub}

\sepsub

\Ssect{Beilinson spectral sequence}{beil}

For every coherent sheaf $\F$ on $\P^2$ there is a free monad, called the
{\em Beilinson free monad}, with middle cohomology $\F$:
\bdm
0 \lra \CC^{-2} \lra \CC^{-1} \lra \CC^0 \lra \CC^1 \lra \CC^2 \lra 0,
\edm
\bdm
\CC^i = \bigoplus_{0 \le j \le n} \H^{i+j}(\F \tensor \Om^j(j)) \tensor \O(-j).
\edm
All maps \ $\H^{i+j}(\F \tensor \Om^j(j)) \tensor \O(-j) \lra \H^{i+j+1}(\F
\tensor \Om^j(j)) \tensor \O(-j)$ \ in the monad are zero. For sheaves with
support of dimension one the Beilinson free monad takes the form

\setcounter{subsub}{1}
\begin{equation}\label{2.E}
0 \lra \CC^{-2} \lra \CC^{-1} \lra \CC^0 \lra \CC^1 \lra 0,
\end{equation}

\begin{align*}
\CC^{-2} & = \H^0(\F(-1)) \tensor \O(-2), \\
\CC^{-1} & = \big(\H^0(\F \tensor \Om^1(1)) \tensor \O(-1) \big) \oplus \big(
\H^1(\F(-1))\tensor \O(-2)\big) , \\
\CC^0 & = \big(\H^0(\F) \tensor \O \big) \oplus\big(\H^1(\F \tensor \Om^1(1))
\tensor\O(-1)\big) , \\
\CC^1 & = \H^1(\F) \tensor \O.
\end{align*}
Dualizing {(\ref{2.E})} we get a free monad for $\F^D$:

\setcounter{subsub}{2}
\begin{equation}\label{2.F}
0 \lra \CC^{-2}_D \lra \CC^{-1}_D \lra \CC^0_D \lra \CC^1_D \lra 0,
\end{equation}
\bdm
\CC^i_D = {\mathcal Hom}(\CC^{-1-i},\omega_{\P^2}).
\edm
For every coherent sheaf $\F$ on $\P^2$ there is a spectral sequence of sheaves,
called the {\em Beilinson spectral sequence}, which converges to $\F$ in degree
zero and to 0 in degree nonzero. Its first term, $\EE^1(\F)$, is given by
\bdm
\EE^1_{ij} = \H^j (\F \tensor \Om^{-i}(-i)) \tensor \O(i).
\edm
If $\F$ is supported on a curve, the relevant part of $\EE^1(\F)$ is exhibited
in the following tableau:
\setcounter{subsub}{3}
\begin{equation}\label{2.G}
\end{equation}
$$\xymatrix@R=12pt{
\H^1(\F(-1))\tensor\O(-2)\fleq[d] & \H^1(\F\tensor\Om^1(1))\tensor\O(-1)\fleq[d]
& \H^1(\F)\tensor\O\fleq[d]\\
\EE_{-2,1}^1\ar[r]^-{\f_1} & \EE_{-1,1}^1\ar[r]^-{\f_2} & \EE_{01}^1
}$$
$$\xymatrix@R=12pt{
\EE_{-2,0}^1\ar[r]^-{\f_3}\fleq[d] & \EE_{-1,0}^1\ar[r]^-{\f_4}\fleq[d] &
\EE_{00}^1\fleq[d]\\
\H^0(\F(-1))\tensor\O(-2) & \H^0(\F\tensor\Om^1(1))\tensor\O(-1)
& \H^0(\F)\tensor\O
}
$$

\SEP

All the other $\EE^1_{ij}$ are zero. The relevant part of $\EE^2$
is
\bdm
\xymatrix
{
\EE^2_{-2,1} = \Ker(\f_1) \ar[rrd]^{\f_5} & \EE^2_{-1,1} = \Ker(\f_2)/\Im(\f_1)
& \EE^2_{01} =\Coker(\f_2) \\
\EE^2_{-2,0} = \Ker(\f_3) & \EE^2_{-1,0} = \Ker(\f_4)/\Im(\f_3) & \EE^2_{00} =
\Coker(\f_4)
}.
\edm
All the other $\EE^2_{ij}$ are zero. The relevant part of $\EE^3$ is
\bdm
\xymatrix
{
\EE^3_{-2,1} = \Ker(\f_5) & \EE^3_{-1,1} = \Ker(\f_2)/\Im(\f_1) & \EE^3_{01}
=\Coker(\f_2) \\
\EE^3_{-2,0} = \Ker(\f_3) & \EE^3_{-1,0} = \Ker(\f_4)/\Im(\f_3) & \EE^3_{00}  =
\Coker(\f_5)
}.
\edm
All the maps in $\EE^3(\F)$ are zero. This shows that $\EE^3=\EE^{\infty}$,
hence all the terms in $\EE^3$, except, possibly, $\EE^3_{00}$ and
$\EE^3_{-1,1}$, are zero. Moreover, there is an exact sequence
\bdm
0 \lra \Coker(\f_5) \lra \F \lra \Ker(\f_2)/\Im(\f_1) \lra 0.
\edm
We conclude that $\f_2$ is surjective and that there are exact sequences

\setcounter{subsub}{4}
\begin{equation}\label{2.H}
0 \lra \H^0(\F(-1)) \tensor \O(-2) \stackrel{\f_3}{\lra} \H^0 (\F \tensor
\Om^1(1)) \tensor \O(-1) \stackrel{\f_4}{\lra} \H^0(\F) \tensor \O
\end{equation}
$$\lra \Coker(\f_4) \lra 0,$$
\setcounter{subsub}{5}
\begin{equation}\label{2.I}
0 \lra \Ker(\f_1) \stackrel{\f_5}{\lra} \Coker(\f_4) \lra \F \lra
\Ker(\f_2)/\Im(\f_1) \lra 0.
\end{equation}

\sEP

Let $S$ be a scheme over $\C$. For every coherent sheaf $\F$ on $\P^2 \times S$
there is a spectral sequence of sheaves on $\P^2 \times S$, called the relative 
Beilinson spectral sequence, which converges to $\F$ in degree zero and to 0 in
degree nonzero. Its $\EE^1$-term is given by
\bdm
\EE^1_{ij} = \R^j_{p_*}(\F \tensor \Om^{-i}(-i)) \boxtimes \O_{\P^2}(i).
\edm
Here \ $p: \P^2 \times S \to S$ \ is the projection onto the second component.

\SEP

\begin{subsub}\label{2.J}{\bf Proposition : }
Let $s$ be a closed point of $S$. If all the base change homomorphisms
\bdm
\R^j_{p_*}(\F \tensor \Om^{-i}(-i))_s \lra \H^j (\F_s \tensor \Om^{-i}(-i))
\edm
are isomorphisms, then the restriction of the $\text{\emph{E}}^1$-term of
the relative Beilinson spectral sequence for $\F$  to \ $\P^2 \times \{ s\}$ \
is the $\text{\emph{E}}^1$-term for the Beilinson spectral sequence for the
restriction $\F_s$ of $\F$ to \ $\P^2 \times \{ s \}$.
\end{subsub}

\Proof  Let $p_1, p_2: \P^2 \times \P^2 \to \P^2$ be the projections
onto the first and second component. We consider the resolution of the diagonal
$\Delta
\subset \P^2 \times \P^2$ given on p. 242 in \cite{oss}:
\bdm
0 \lra \O_{\P^2}(-2) \boxtimes \Om_{\P^2}^2 (2) \lra \O_{\P^2}(-1) \boxtimes
\Om_{\P^2}^1(1) \lra \O_{\P^2 \times \P^2} \lra \O_{\Delta} \lra 0.
\edm
The maps in $\EE^1(\F)$ are the induced maps
\begin{align*}
\EE^1_{ij}(\F) =
& \R^j_{(p_1 \times 1_S)_*}\big(\O_{\P^2}(i) \boxtimes (p_2 \times 1_S)^*(\F
\tensor \Om_{\P^2}^{-i}(-i))\big) \lra  \\
& \R^j_{(p_1 \times 1_S)_*}\big(\O_{\P^2}(i+1) \boxtimes (p_2 \times 1_S)^*(\F
\tensor \Om_{\P^2}^{-i-1}(-i-1))\big) = \EE^1_{i+1,j}(\F).
\end{align*}
Restricting to $\P^2 \times \{ s \}$ we get the map
\bdm
\O_{\P^2}(i) \tensor \R^j_{p_*}(\F \tensor \Om_{\P^2}^{-i}(-i))_s \lra
\O_{\P^2}(i+1) \tensor\R^j_{p_*}(\F \tensor \Om_{\P^2}^{-i-1}(-i-1))_s.
\edm
From the naturality of the base-change homomorphism we see that the above is the
induced map
\bdm
\O_{\P^2}(i) \tensor \H^j (\F_s \tensor \Om_{\P^2}^{-i}(-i)) \lra \O_{\P^2}(i+1)
\tensor\H^j (\F_s \tensor \Om_{\P^2}^{-i-1}(-i-1)).
\edm
But this is the map from $\EE^1(\F_s)$, which finishes the proof of the claim.
\carre

\SEP

\begin{subsub}\label{2.K}{\bf Proposition : }
Let $S$ be a noetherian integral scheme over $\C$ and let
$\F$ be a coherent sheaf on $\P^2 \times S$ which is $S$-flat. For a closed
point $s$ in $S$ we denote by $\F_s$ the restriction of $\F$ to $\P^2 \times \{
s \}$. Assume that for all $i$ and $j$, $h^j(\F_s \times \Om^{-i}(-i))$ is
independent of $s$. Then, for all closed points $s$ in $S$, the restriction of
$\text{\emph{E}}^1(\F)$ to $\P^2 \times \{ s\}$ is
$\text{\emph{E}}^1(\F_s)$.
\end{subsub}

\Proof  According to III 12.9 from \cite{hartshorne}, all base change
homomorphisms from \ref{2.J} are isomorphisms, so \ref{2.K} is a corollary of
\ref{2.J}.
\carre
\end{sub}

\sepsub

\Ssect{Quotients by reductive groups and moduli spaces of sheaves}{mod}

We first recall the definition of good and geometric quotients (cf.
\cite{mumf}, \cite{newstead}) :

\SEP

\begin{subsub}\label{2.L0}{\bf Definition : }
Let an algebraic group $G$ act on an algebraic variety $X$. Then a pair
$(\varphi, Y)$ of a variety and a morphism $X\xrightarrow{\varphi} Y$ is called
a good quotient if
\begin{enumerate}
\item [(i)] $\varphi$ is $G$--invariant (for the trivial action of $G$ on
  $Y$),
\item [(ii)] $\varphi$ is affine and surjective,
\item [(iii)] If $U$ is an open subset of $Y$ then $\varphi^\ast$ induces an
  isomorphism $\O_Y(U)\simeq\O_X(\varphi^{-1}U)^G$, where the latter
  denotes the ring of $G$--invariant functions,
\item [(iv)] If $F_1, F_2$ are disjoint closed and $G$--invariant subvarieties
  of $X$ then $\varphi(F_1)$, $\varphi(F_2)$ are closed and disjoint.
\end{enumerate}

If in addition the fibres of $\varphi$ are the orbits of the action the quotient
$(\varphi, Y)$ is called a geometric quotient.
\end{subsub}

\SEP

C.~Simpson's construction of the moduli spaces $\M(r,\chi)$ (cf.
\cite{lepotier}, \cite{si}) is based on the following facts: there are a smooth
variety $R$ and a reductive group $G$ (to be precise, $G$ is a special linear
group) acting algebraically on $R$, such that $\M(r,\chi)$ is a good quotient of
$R$ by $G$. Moreover, the open subset $\M^{\text{s}}(r,\chi)$ of isomorphism
classes of stable sheaves is the geometric quotient of an open subset $R_0
\subset R$ modulo $G$. There is a coherent $R$-flat sheaf $\tilda{\F}$ on $\P^2
\times R$ whose restriction $\tilda{\F}_s$ to every closed fiber $\P^2 \times \{
s\}$ is a semistable sheaf with Hilbert polynomial $P(t)=rt+\chi$. The quotient
morphism
$$\pi : R \lra \M(r,\chi)$$
maps $s$ to the stable-equivalence class of $\tilda{\F}_s$, denoted
$[\tilda{\F}_s]$. The quotient morphism \mM{R_0 \to \M^{\text{s}}(r,\chi)} sends
$s$ to the isomorphism class of $\tilda{\F}_s$.

\SEP

\begin{subsub}\label{2.L}{\bf Proposition : }
Let $X$ be an irreducible locally closed subvariety
of $\M(r,\chi)$ and $S' \subset R$ the preimage of $X$ equipped with the
canonical reduced induced structure. Then

(i) The restriction of $\pi:S'\to X$ is a good quotient of $S'$ by $G$.

(ii) There exists an irreducible component $S$ of $S'$ such that $\pi(S)=X$.

(iii) The restriction of $\pi:S\to X$ is a good quotient of $S$ by $G$.
\end{subsub}

\Proof 
It follows from \cite{sw}, 3.2 (i), that $\pi^{-1}(X)\to X$ is a good quotient
by $G$. We have \mM{S'=\pi^{-1}(X)_{red}}. Since $\pi$ is affine (i) is
reduced to the following~: let $Z={spec}(A)$ be an affine scheme and suppose
given an algebraic action of $G$ on $X$. Suppose that $Z//G={spec}(A^G)$ is
reduced. Then the canonical morphism \hbox{$\phi:A^G\to(A/\text{rad}(A))^G$} is
an
isomorphism. We have $\Ker(\phi)=A^G\cap\text{rad}(A)$, and since $A^G$ is
reduced, $\phi$ is injective. We have a commutative square
\[
\xymatrix{A\flon[d]\flon[r]^R & A^G\ar[d]^\phi \\
A/\text{rad}(A)\flon[r]^R & (A/\text{rad}(A))^G}\]

\noindent
$R$ beeing the Reynolds operators (cf. \cite{mumf}), showing that $\phi$ is
surjective. This proves (i).

Let $S''$ be an irreducible component of $X$. We have \ $S''\subset G.S''$,
and  $G.S''$ is irreducible (because $G$ and $S''$ are). It follows
that $GS''=S''$ and that $\pi(S'')$ is closed in $X$. Since $X$ is irreducible
and the union of the closed images of the irreducible components of $S'$ there
exists an irreducible component $S$ of $S'$ such that $\pi(S)=X$. This proves
(ii).

Again in the case of (iii) the problem is local. So we can suppose that
$S'={spec}(A)$ and $X={spec}(A^G)$. Let $I$ be the $G$-invariant ideal of $S$.
Then we have to prove that the canonical morphism $\phi:A^G\to(A/I)^G$ is an
isomorphism. The surjectivity of $\pi:S\to X$ implies that the composition
$$A^G \stackrel{\phi}{\lra} (A/I)^G\subset A/I$$
is injective, so $\phi$ is injective. The surjectivity of $\phi$ can be seen
using the Reynolds operators as before. This proves (iii).
\carre

\SEP

\begin{subsub}\label{2.M}{\bf Proposition : }
With the above notations, let $\tilda{\F}_{S}$ be the
restriction
of $\tilda{\F}$ from $\P^2 \times R$ to $\P^2 \times S$. Assume that for all $i$
and $j$, $h^j(\tilda{\F}_s \tensor \Om^{-i}(-i))$ is independent of the closed
point $s$ in $S$.
Then, for all closed points $s$ in $S$, the restriction of
$\text{\emph{E}}^1(\tilda{\F})$ to $\P^2 \times \{ s\}$ is
$\text{\emph{E}}^1(\tilda{\F}_s)$.
\end{subsub}

\Proof The hypotheses of \ref{2.K} are satisfied because flatness is
preserved under pulling back. \carre
\end{sub}

\sepsub

\Ssect{Kronecker modules}{KM}

Let $L$ be a finite dimensional nonzero vector space over $\C$, and $m$, $n$
positive integers. Let $W=L(\C^m\otimes L,\C^n)\backslash\lbrace 0\rbrace$. We
have an algebraic action of
$$\Gamma=(\text{GL}(m)\times\text{GL}(n)/\C^*$$
on $W$ given by
$$\xymatrix@R=4pt{\Gamma\times W\ar[r] & W\\
((g_1,g_2),f)\fmaps[r] & g_2\circ f\circ(g_1\otimes I_L)^{-1}
}$$
Let $\P=\P(L(\C^m\otimes L,\C^n))$. The preceeding action induces an action of
the reductive group \ $G=\text{SL}(m)\times\text{SL}(m)$ \ on $\P$ with an
obvious linearization. According to K.~Hulek \cite{hu}, the linear maps
$\C^m\otimes L\to\C^n$ are called {\em $L$-Kronecker modules}. A Kronecker
module will be called {\em semi-stable} (resp. {\em stable}) if it is nonzero
and if the corresponding point in $\P$ is semi-stable (resp. stable) for the
above action. We have

\SEP

\begin{subsub}\label{2.N}{\bf Proposition : }
A $L$-Kronecker module \ $\tau:\C^m\otimes L\to\C^n$ is semi-stable (resp.
stable) if and only for every linear subspaces $H\subset\C^m$, $K\subset\C^n$,
with $H\not=\lbrace 0\rbrace$, such that \ $\tau(H\otimes L)\subset K$ \ we have
\[\frac{\dim(K)}{\dim(H)} \ \geq \ \frac{n}{m} \quad\quad
\text{(resp.} \quad > \text{)} .\]
\end{subsub}

(cf. \cite{dr2}, prop. 15, \cite{ki}).

\SEP

Let $\P^{ss}$ (resp. $\P^{s}$) denote the $G$-invariant open subset of
semi-stable (resp. stable) points of $\P$. Let 
\[N(L,m,n)=\P^{ss}//G , \quad N_s(L,m,n)=\P^{s}/G .\]
Of course these varieties depend only on $m$, $n$ and $\dim(L)$.
If $\dim(L)=q$ we will also use the notations $N(q,m,n)$ and $N_s(q,m,n)$.

The variety $N(q,m,n)$ is projective, irreducible and locally factorial, and
$N_s(q,m,n)$ is a smooth open subset of $N(q,m,n)$.

Let $x_q$ be the smallest solution of the equation \ $X^2-qX+1=0$ . Then we
have \hfil\break $\dim(N(q,m,n))>0$ \ if and only if \
$x_q<\frac{m}{n}<\frac{1}{x_q}$ .
In this case $N_s(q,m,n)$ is not empty and we have \
$\dim(N(q,m,n))=qmn-m^2-n^2+1$ .

If $m$ and $n$ are relatively prime then \ $N(q,m,n)=N_s(q,m,n)$ \ hence
$N(q,m,n)$ is a projective smooth variety. In this case there exists a {\em
universal morphism} on $N(q,m,n)$:  there are algebraic vector bundles $E$,
$F$ on $N(q,m,n)$ of rank $m$, $n$ respectively, and a morphism \ $\tau:E\otimes
L\to F$ \ such that for every closed point $x$ of $N(q,m,n)$, and isomorphisms
$E_x\simeq\C^m$, $F_x\simeq\C^n$, the linear map \ $\tau_x:\C^m\otimes
L\to\C^n$ \ belongs to the $G$-orbit represented by $x$.

The moduli spaces of Kronecker modules appear in the following context :
suppose given two vector bundles $U$, $V$ on $\P^2$. Then a morphism
$U\otimes\C^m\to V\otimes\C^n$ is equivalent to a $\text{Hom}(U,V)^*$-Kronecker
module $\text{Hom}(U,V)^*\otimes\C^m\to \C^n$.
\end{sub}

\sepsub

\Ssect{Moduli spaces of morphisms}{MSM}

Let $X$ be a projective algebraic variety, and $r$, $s$ positive integers. For
\mM{1\leq i\leq r} (resp. \mM{1\leq j\leq s}) let $m_i$ (resp. $n_j$) be a
positive integer and $\F_i$ (resp. $\E_j$) a coherent sheaf on $X$. Let
\[
\E=\underset{1\leq i\leq r}{\bigoplus}\E_i\otimes\C^{m_i}\quad\text{
and }\quad\F = \underset{1\leq j\leq s}{\bigoplus} \F_j\otimes\C^{n_j} .
\]
We suppose that the sheaves $\E_i$, $\F_j$ are simple and that
\[\Hom(\E_i,\E_k)=\Hom(\F_j,\F_l)=\nsp\]
if $i>k$ and $j>l$.

Let $\W=\Hom(\E,\F)$. The algebraic group \ $G=\Aut(\E)\times\Aut(\F)$ \ acts on
$\W$ in an obvious way. We can see the elements of $\Aut(\E)$ as matrices
\[
\left (
\begin{array}{ccccc}
g_1 & 0 & \cdots & 0\\
u_{21} & g_2 & & \vdots\\
\vdots & \ddots & \ddots & 0\\
u_{r1} & \cdots & u_{r, r-1} & g_r
\end{array}
\right)
\]
where $g_i\in \GL(m_i)$ and \ $u_{ki} \in L(\C^{m_i},\C^{m_k})\otimes \Hom(\E_i,
\F_k)$ \ (and similarly for $\Aut(\F)$).
Let \ $G_{red}=\prod\ \GL(m_i)\times\prod\GL(n_l)$ , which is a reductive
subgroup of $G$, and $H$ the maximal normal unipotent subgroup of $G$,
consisting of pairs of matrices with identities as diagonal terms.

The action of $G_{red}$ is well known (cf. \cite{ki}). Let
\ $\sigma=(\lambda_1,\ldots,\lambda_r,\mu_1,\ldots,\mu_s)$  \ be a
sequence of positive rational numbers such that
\[\sigg_{1\leq i\leq r}\lambda_i m_i \ = \ \sigg_{1\leq j\leq s}\mu_j n_j
\ = \ 1 \]
(such a sequence is called a {\em polarization}).
An element $f\in\W$ is called {\em $G_{red}$--semi-stable} (resp. {\em
$G_{red}$--stable}) with respect to $\sigma$ if for any choice of subspaces
$M_i\subset\C^{m_i}$, $N_j\subset\C^{n_j}$ such that $N_j\not=\C^{m_j}$ for at
least one $j$, and that $f$ maps $\oplus(\E_i\otimes M_i)$ into
$\oplus(\F_j\otimes N_j)$, we have
\[\sigg_{1\leq i\leq r}\lambda_i\dim(M_i)\leq\sigg_{1\leq j\leq s}\mu_j\dim(N_j)
\quad \text{(resp. }<\text{)} . \]
There exists a good quotient of the open subset of $G_{red}$--semi-stable points
of $\W$ with respect to $\sigma$.

We consider now the action of the whole group $G$ which is not reductive
in general. An element $f\in\W$ is called {\em $G$--semi-stable} (resp. {\em
$G$--stable}) with respect to $\sigma$ if all the elements of $H.f$ are
$G_{red}$--semi-stable (resp. $G_{red}$--stable) with respect to $\sigma$. Let
$\W^{ss}(\sigma)$ (resp. $\W^s(\sigma)$) be the open $G$-invariant subset of
$G$--semi-stable (resp. $G$--stable) points of $\W$ with respect to $\sigma$.
If suitable numerical conditions are satisfied by $\sigma$ then
$\W^{ss}(\sigma)$ admits a good and projective quotient and $\W^s(\sigma)$
admits a geometric quotient, which is smooth (cf. \cite{drezet-trautmann}).

\end{sub}

\sepsub


\sepsec

\section{Euler Characteristics One and Three}\label{M4143}

\Ssect{The open strata}{open1}

According to 4.2 in \cite{maican}, the sheaves $\G$ giving a point in $\M(4,3)$
and satisfying \hfil\break\mM{h^0(\G(-1))=0} are precisely the sheaves that have
a resolution of the form
\setcounter{subsub}{1}
\begin{equation}\label{3.A1}
0 \lra \O(-2) \oplus 2\O(-1) \stackrel{\f}{\lra} 3\O \lra \G \lra 0
\end{equation}
with $\f_{12}$ having all maximal minors nonzero.

\SEP

\begin{subsub}\label{3.A0} {\bf Moduli spaces of morphisms -} \rm
It is easy to see, using the stability conditions of \ref{2.N}, that a morphism
$f:2\O(-1)\to 3\O$ is stable (as a $V$-Kronecker module) if and only if all its
maximal minors are nonzero. Moreover a stable morphism is injective (as a
morphism of sheaves).

We consider now morphisms as in (3.1.1)
\[\f : \O(-2)\oplus 2\O(-1)\lra 3\O .\]
Let $W=\text{Hom}(\O(-2)\oplus 2\O(-1),3\O)$. Then the linear algebraic group
$$G=\big(\Aut(\O(-2)\oplus 2\O(-1))\times\Aut(3\O)\big)/\C^*$$
acts on $W$ in an obvious way. Good quotients of some $G$-invariant open
subsets of $W$ are given in \cite{drezet-trautmann}, 9.3 (cf. \ref{MSM}). The
quotient related to $\M(4,3)$ is the obvious one, and we will describe it.

We begin with the following remark : let \ $f:2\O(-1)\to 3\O$ \ be an injective
morphism, and $\sigma\in \H^0(\Coker(f)(2))$. Then $\sigma$ can be lifted to a
morphism $\O(-2)\to 3\O$ and defines thus with $f$ a morphism as in
(3.1.1). All the morphisms constructed in this way are in the same $G$-orbit,
and every morphism $\f:\O(-2)\oplus 2\O(-1)\to 3\O$ \ such that
$\f_{12}=f$ comes from a section of $\Coker(f)(2)$.

Let $\tau:E\times V\to
F$ be a universal morphism on $N(3,2,3)$ (cf. \ref{KM}). Let $p_1$, $p_2$ be the
projections \ $N(3,2,3)\times\P^2\to N(3,2,3)$, $N(3,2,3)\times\P^2\to\P^2$
respectively. From $\tau$ we get a morphism of sheaves on \ $N(3,2,3)\times\P^2$
$$\theta : p_1^*(E)\otimes p_2^*(\O(-1))\otimes V\lra p_1^*(F)
\ .$$
For every $x\in N(3,2,3)$, this morphism is injective on the fiber \ $\lbrace
x\rbrace\times\P^2$. Hence $\Coker(\theta)$ is flat on $N(3,2,3)$, and for
every $x\in N(3,2,3)$ we have $\Coker(\theta)_x=\Coker(\theta_x)$.
Let
$$U \ = \ p_{1*}(\Coker(\theta)\otimes p_2^*(\O(2))) .$$
It is a rank 3 vector bundle on $N(3,2,3)$, and for every $x\in N(3,2,3)$ we
have \hfil\break $U_x=\H^0(\Coker(\theta_x)(2))$.

Let $\mathbb{W}=\P(U)$. Let $\boldsymbol W$ be the open $G$-invariant subset of
$W$ consisting of morphisms $\f$ such that $\f_{12}$ is stable, and such that
the section of $\Coker(\f_{12})$ defined by $\f_{22}$ is nonzero. We have
an obvious morphism $\boldsymbol W\to \mathbb{W}$ which is a geometric
quotient, i.e. $\mathbb{W}=\boldsymbol W/G$. Hence $\boldsymbol W/G$ is a
smooth projective variety.
\end{subsub}

\SEP

\begin{subsub}\label{3.A2} {\bf The open stratum of $\M(4,3)$ -} \rm Let
$X_0(4,3)$ be the open subset of $\M(4,3)$ corresponding to sheaves $\G$ such
that $h^0(\G(-1))=0$. Let $W_0(4,3)\subset\boldsymbol W$ be the $G$-invariant
open subset consisting of injective morphisms. We will see later that the
morphism $\rho_{4,3}:W_0(4,3)\to X_0(4,3)$ sending a morphism to its cokernel is
a geometric quotient by $G$, and hence $X_0(4,3)$ is isomorphic to the open
subset of $\mathbb{W}$ corresponding to injective morphisms.
\end{subsub}

\SEP

Now we will describe the open stratum of $\M(4,1)$. As at 2.2.2, dualizing the
exact sequence (3.1.1) we get a resolution for the sheaf $\F=\G^D(1)$. According
to \ref{2.A}, the latter is in $\M(4,1)$. Hence we obtain

\SEP

\begin{subsub}\label{3.A}{\bf Proposition : }
The sheaves $\F$ in $\MM(4,1)$ satisfying $h^1(\F)=0$
are precisely the sheaves with resolution of the form
\bdm
0 \lra 3\O(-2) \stackrel{\f}{\lra} 2\O(-1) \oplus \O \lra \F \lra 0,
\edm
where all maximal minors of $\f_{11}$ are nonzero.
\end{subsub}

\SEP

Let $W_0=W_0(4,1)$ denote the set of injective morphisms
\bdm
\f : 3\O(-2) \lra 2\O(-1) \oplus \O,
\edm
such that $\f_{11}$ is stable. The linear algebraic group
$$\big(\Aut(3\O(-2)) \times \Aut(2\O(-1)\oplus\O)\big)/\C^*$$
acts on $\text{Hom}(3\O(-2),2\O(-1)\oplus\O)$ by conjugation. Of course this
group is canonically isomorphic to $G$, and the isomorphism
\[\xymatrix@R=3pt{\text{Hom}(\O(-2)\oplus 2\O(-1),3\O)\ar[r]^\Lambda &
\text{Hom}(3\O(-2),2\O(-1)\oplus\O)\\
\f\fmaps[r] & {}^t\f\otimes I_{\O(-2)}}\]
is $G$-invariant and sends $W_0(4,3)$ to $W_0(4,1)$.

\SEP

\begin{subsub}\label{3.A3}{\bf The open stratum of $\M(4,1)$ -} \rm
Let $X_0=X_0(4,1)$ be the open subset of $\M(4,1)$ of sheaves $\F$ such that
$h^1(\F)=0$. Let
\[\xymatrix@R=3pt{\lambda : \M(4,3)\ar[r]^\simeq & \M(4,1)\\
\G\fmaps[r] & \G^D(1)}\]
be the isomorphism of \ref{2.A}.
\end{subsub}

The map \ $\rho : W_0 \to X_0$ \ (resp. \ $\rho_{4,3}:W_0(4,3)\to
X_0(4,3)$) which sends $\f$ to the isomorphism class
of $\Coker(\f)$ is a surjective morphism whose fibers are $G$-orbits. We have a
commutative diagram
\[\xymatrix{W_0(4,3)\ar[r]^\Lambda\ar[d]^{\rho_{4,3}} & W_0(4,1)\ar[d]^\rho\\
X_0(4,3)\ar[r]^\lambda & X_0}\]

We claim that $\rho$ (anf hence also $\rho_{4,3}$) is a geometric quotient map:

\SEP

\begin{subsub}\label{3.B}{\bf Theorem : }
The geometric quotient $W_0/G$ is isomorphic to $X_0$.
\end{subsub}

\Proof We will show that $\rho$ is a categorical quotient map and the
isomorphism \hfil\break$W_0/G \isom X_0$ will follow from the uniqueness of the
categorical quotient. Given a $G$-invariant morphism of varieties $f: W_0 \to
Y$, there is a unique map $g: X_0 \to Y$ such that $g \circ \rho = f$. We need
to show that $g$ is a morphism of varieties. To see this we consider the good
quotient $\pi: S \to X_0$ of \ref{2.L}. We will show that $S$ can
be covered with open sets $U$ for which there are morphisms $\s_U : U \to W_0$
making the diagram commute:
\bdm
\xymatrix
{
U \ar[rr]^{\s_U} \ar[rd]_{\pi} & & W_0 \ar[ld]^{\rho} \\
& X_0
}.
\edm
Now we note that \ $g \circ \pi : S \to Y$ is a morphism because its restriction
to each open set $U$ is $f \circ \s_U$. Thus $g$ is the unique morphism
associated to $g \circ \pi$ by the categorical quotient property of $\pi$.

It remains to construct the morphisms $\s_U$. This can be achieved as follows:
we notice that $\tilda{\F}_S$ satisfies the hypotheses of \ref{2.M}.
This is so because for every sheaf $\F$ giving a point in $X_0$ we have
\begin{align*}
& h^1(\F \tensor \Om^2(2))=3, \qquad \qquad && h^1(\F \tensor \Om^1(1))=2,
\qquad \qquad &&& h^1(\F)=0, \\
& h^0(\F \tensor \Om^2(2))=0, \qquad \qquad && h^0(\F \tensor \Om^1(1))=0,
\qquad \qquad &&& h^0(\F)=1.
\end{align*}
Thus the higher direct images $\R^j_{p_*}(\tilda{\F}_S \tensor \Om^{-i}(-i))$
are locally free sheaves on $S$; we cover $S$ with open subsets $U$ on which
they are trivial and we fix such trivializations. For an arbitrary closed point
$s$ in $U$ we restrict $\EE^1(\tilda{\F}_S)$ to $\P^2 \times \{ s \}$.
In this manner we construct a morphism $\z_U$ from $U$ to an open subset
$E$ of the space of spectral sequences with $\EE^1$-term
\bdm
\xymatrix
{
\EE^1_{-2,1}= 3\O(-2) \ar[r]^{\f_1} & \EE^1_{-1,1}= 2\O(-1) \ar[r]^{\f_2} &
\EE^1_{01} =0 \\ \EE^1_{-2,0}=0 & \EE^1_{-1,0}=0 & \EE^1_{00}=\O
}.
\edm
All the other $\EE^1_{ij}$ are zero. By \ref{2.M} $\z_U (s)$ is isomorphic to
$\EE^1(\tilda{\F}_s)$.  It remains to construct a morphism $\xi : E \to W_0$
which maps $\z_U(s)$ to a point in $\f \in W_0$ satisfying $\Coker(\f) \isom
\tilda{\F}_s$; then we can put $\s_U = \xi \circ \z_U$. In other words, we need
to obtain a resolution of the form \ref{3.A} for $\tilda{\F}_s$ starting
with
$\EE^1(\tilda{\F}_s)$ and performing algebraic operations.

For the problem at hand it is easier to construct $\xi$ for the dual sheaf
$\G=\tilda{\F}_s^D(1)$. In view of \ref{2.A} and (2.2.2) we can dualize the
problem: given $\G$ in $\M(4,3)$ with $h^0(\G(-1))=0$, we would like to
construct the dual to resolution \ref{3.A} starting from $\EE^1(\G)$.
Tableau (2.2.3) is
\bdm
\xymatrix
{
\O(-2) & 0 & 0 \\
0 & 2\O(-1) \ar[r]^{\f_4} & 3\O
}
\edm
and the exact sequence (2.2.5) takes the form
\bdm
0 \lra \O(-2) \stackrel{\f_5}{\lra} \Coker(\f_4) \lra \G \lra 0.
\edm
It is now clear that $\G$ is isomorphic to the cokernel of the morphism
\bdm
\f = (\f_5,\f_4) : \O(-2) \oplus 2\O(-1) \lra 3\O.
\edm
This finishes the proof of \ref{3.B}. \carre
\end{sub}

\sepsub

\Ssect{The closed strata}{clos1}

We will treat only the closed stratum in $\M(4,1)$, the case of $\M(4,3)$ can
be deduced by duality as in \ref{open1}.
According to \ref{2.C}, we have $h^0(\F(-1))=0$ for all $\F$ in $\M(4,1)$.
From this and 5.3 in \cite{maican} we obtain:

\SEP

\begin{subsub}\label{3.C}{\bf Theorem : }
The sheaves $\F$ in $\MM(4,1)$ satisfying $h^1(\F)=1$
are precisely the sheaves with resolution of the form
\bdm
0 \lra \O(-3) \oplus \O(-1) \stackrel{\f}{\lra} 2\O \lra \F \lra 0,
\edm
where $\f_{12}$ and $\f_{22}$ are linearly independent vectors in $V^*$.
\end{subsub}

\SEP

The description of the closed strata is very similar to that of
the open ones. Let \hfil\break\mM{W_1=W_1(4,1)} denote the set of morphisms from
\ref{3.C} that is, the set of injective morphisms
\bdm
\f : \O(-3) \oplus \O(-1) \lra 2\O,
\edm
such that $\f_{12}$ and $\f_{22}$ are linearly independent.
The linear algebraic group
$$G=\big(\Aut(\O(-3) \oplus \O(-1)) \times \Aut(2\O)\big)/\C^*$$
acts on $W_1$ by conjugation. Note that a morphism $f:\O(-1)\to 2\O$ is stable
(as a $V$-Kronecker module) if and only if $f_{12}$ and $f_{22}$ are linearly
independent, and a stable morphism is injective. The variety $N(3,1,2)$ is
canonically isomorphic to $\P^2$: to every point $P$ of $\P^2$, represented by
a line $D\subset V$, the corresponding stable morphism is the canonical one
$\O(-1)\to\O\otimes (V/D)$. The cokernel of this morphism is
isomorphic to $\I_P(1)$ ($\I_P$ beeing the ideal sheaf of $P$). 

As for the open strata, given a stable $f:\O(-1)\to 2\O$ with cokernel
$\I_P(1)$, the morphisms $\O(-3)\oplus\O(-1)\to 2\O$ whose restriction to
$\O(-1)$ is $f$ correspond to sections of $\I_P(4)$. The morphisms
corresponding to nonzero sections are injective. Let $Q_4$ denote the cokernel
of the canonical morphism $\O(-4)\to\O\otimes S^4V$ (cf. \cite{dr2})
which is a rank 14 vector bundle on $N(3,1,2)=\P^2$. We have a canonical
isomorphism $Q_{4P}^*\simeq\H^0(\I_P(4))$, for every point $P$ of $\P^2$.
Let $\mathbb{W}'=\P(Q_4^*)$. Then it is easy to see that we obtain a geometric
quotient \
$W_1\to\mathbb{W}'$.

\SEP

Let $X_1=X_1(4,1)$ be the locally closed subset of $\M(4,1)$ given by the
condition \hfil\break$h^1(\F)=1$. We equip $X_1$ with the canonical induced
reduced structure.

\SEP

\begin{subsub}\label{3.D}{\bf Theorem : }
The geometric quotient $\mathbb{W}'$ is isomorphic to $X_1$.
In particular, $X_1$ is a smooth closed subvariety of codimension 2.
\end{subsub}

\Proof As explained in the proof of \ref{3.B}, we have to construct the
resolution of \ref{3.C} starting from $\EE^1(\F)$. Tableau (2.2.3) is
\bdm
\xymatrix
{
3\O(-2) \ar[r]^{\f_1} & 3\O(-1) \ar[r]^{\f_2} & \O \\
0 & \O(-1) \ar[r]^{\f_4} & 2\O
}.
\edm
As $\f_2$ is surjective, we may assume that it is given by the matrix
${\displaystyle
\begin{pmatrix}
X & Y & Z
\end{pmatrix}
}$.
Thus $\f_1$ has columns that are linear combinations of the columns of the
matrix
\bdm
\begin{pmatrix}
-Y & -Z & 0 \ \\
X & 0 \ & -Z \\
0 \ & X & Y
\end{pmatrix}
.
\edm
As $\F$ surjects onto $\Ker(\f_2)/\Im(\f_1)$, the latter has rank zero, hence
$\Im(\f_1)$ has rank two, hence the columns of $\f_1$ span a two or three
dimensional vector space. In the former case we have the isomorphism
$\Ker(\f_2)/\Im(\f_1) \isom\Om^1/2\O(-2)$. This sheaf has Hilbert polynomial
$P(t)=t-2$, which contradicts the semistabilty of $\F$. We conclude that $\f_1$
has three linearly independent columns, hence $\Ker(\f_2) = \Im(\f_1)$,
$\Ker(\f_1) \isom \O(-3)$ and (2.2.5) takes the form
\bdm
0 \lra \O(-3) \stackrel{\f_5}{\lra} \Coker(\f_4) \lra \F \lra 0.
\edm
It follows that $\F$ is isomorphic to the cokernel of the morphism
\bdm
\f=(\f_5,\f_4): \O(-3) \oplus \O(-1) \lra 2\O.
\edm
This finishes the proof of \ref{3.D}. \carre

\SEP

\begin{subsub}\label{3.E}{\bf Theorem : }
There are no sheaves $\F$ in $\MM(4,1)$ satisfying the relation
\hfil\break$h^1(\F)\ge 2$.
\end{subsub}

\Proof Let $\F$ be an arbitrary sheaf in $\M(4,1)$. The Beilinson monad (2.2.1)
for $\F(-1)$ has the form
\bdm
0 \lra 7\O(-2) \lra 10\O(-1) \lra 3\O \lra 0.
\edm
Thus $\M(4,1)$ is parametrized by an open subset $M$ inside the space of monads
of the form
\bdm
0 \lra 7\O(-1) \stackrel{A}{\lra} 10\O \stackrel{B}{\lra} 3\O(1) \lra 0.
\edm
Let \ $\eta: M \to \M(4,1)$ \ be the surjective morphism which sends a monad to
the isomorphism class of its cohomology. The tangent space of $M$ at an
arbitrary point $(A,B)$ is
\bdm
{\mathbb{T}}_{(A,B)}M = \{ (\a,\b), \quad \b \circ A +B \circ \a = 0 \}.
\edm
Let \ $\Phi: M \to \Hom(10\O,3\O(1))$ \ be the projection onto the second
component. We claim that $\Phi$ has surjective differential at every point.
Indeed, $\text{d}\Phi_{(A,B)}(\a,\b)=\b$, so we need to show that, given $\b$,
there is $\a$ such that $\b \circ A + B \circ \a = 0$. This follows from the
surjectivity of the map
\bdm
\Hom(7\O(-1),10\O) \lra \Hom(7\O(-1),3\O(1)), \qquad \a \lra B \circ \a.
\edm
To see this we apply the long $\Ext(7\O(-1),\_\_)$ sequence to the exact
sequence
\bdm
0 \lra \Ker(B) \lra 10\O \stackrel{B}{\lra} 3\O(1) \lra 0
\edm
and we use the vanishment $\Ext^1(7\O(-1),\Ker(B))=0$.
The latter follows from the exact sequence
\bdm
0 \lra 7\O(-1) \lra \Ker(B) \lra \F \lra 0
\edm
and the vanishment $\H^1(\F(1))=0$. The latter follows from \ref{2.C}.

\noindent \\
Note that $h^0(\F)=10-\text{rank}(B)$. The subset $N \subset M$ of monads with
cohomology sheaf $\F$ satisfying $h^1(\F) \ge 2$ is the preimage under $\Phi$
of the set of morphisms of rank at most 7. Since any matrix of rank 7 is the
limit of a sequence of matrices of rank 8, and since the derivative of $\Phi$ is
surjective at every point, we deduce that $N$ is included in
$\overline{\eta^{-1}(X_1)} \setminus \eta^{-1}(X_1)$. But, according to
\ref{3.D}, $X_1$ is closed. We conclude that $N$ is empty, which proves
\ref{3.E}. \carre
\end{sub}

\sepsub

\Ssect{Description of the strata}{desc1}

We will describe the open stratum $X_0(4,3)$ and the closed stratum $X_1(4,1)$.
One can obtain a description of the two other strata by duality.
We consider first the open stratum $X_0(4,3)$.

Let $Y$ denote the closed subvariety of $N(3,2,3)$ corresponding to morphisms
\hfil\break$f:2\O(-1)\to 3\O$ such that $\Coker(f)$ is not torsion free, and
$\widetilde{Y}$ the closed subvariety of $X_0(4,3)$ of points over $Y$.

From \cite{dr3} prop. 4.5, we have

\SEP

\begin{subsub}\label{3.F}{\bf Proposition : }Let $f:2\O(-1)\to 3\O$ be a stable
morphism. Then
$\Coker(f)$ is torsion free if and only if it is isomorphic to $\I_Z(2)$, $Z$
beeing a finite subscheme of length 3 of $\P^2$ not contained in any line, and
$\I_Z$ its ideal sheaf.
\end{subsub}

\SEP

Let $f:2\O(-1)\to 3\O$ be a morphism as in \ref{3.F}, i.e. such that $\Coker(f)$
is torsion free. Let \ $\f:\O(-2)\oplus
2\O(-1)\to 3\O$ \ such that $\f_{12}=f$, and $\F=\Coker(\f)$. Then $\f$ is
injective and we have an exact sequence
\[0\lra\O(-2)\lra\I_Z(2)\lra\F\lra 0 .\]
Let $C$ be the quartic defined by the composition
$\O(-2)\to\I_Z(2)\subset\O(2)$, this quartic contains $Z$. Then from the
preceeding exact sequence we get the following one
\[0\lra\F\lra\O_C(2)\lra\O_Z\lra 0 .\]
It follows easily that

\SEP

\begin{subsub}{\bf Proposition : } The open subset
$X_0(4,3)\backslash\widetilde{Y}$ consists of the
kernels of the surjective morphisms \mM{\O_C(2)\to\O_Z}, where $C$ is a quartic
and $Z$ a length 3 finite subscheme of $\P^2$ not contained in any line.

The generic point of $X_0(4,3)$ is of the form $\O_C(2)(-P-Q-R)$, where $C$ is
a smooth quartic and $P$, $Q$, $R$ are points of $C$ not contained in the same
line.
\end{subsub}

\SEP

\begin{subsub} {\bf Non torsion free cokernels -} \rm Let $f:2\O(-1)\to 3\O$ be
a stable morphism such that $\Coker(f)$ is not torsion free. This is the case if
and only if the maximal minors of $f$ have a common divisor, which is a linear
form. It follows easily that there exists a basis $(z_0,z_1,z_2)$ of $V^*$ such
that $f$ is given (up to the action of $\text{GL}(2)\times\text{GL}(3)$) by the
matrix $\begin{pmatrix}z_1 & -z_0\\ z_2 & 0\\ 0 & z_2\end{pmatrix}$. Here all
the maximal minors are multiples of $z_2$. Let $D\subset V$ be the plane defined
by $z_2$ and $\ell$ the corresponding line of $\P^2$. Then $f$ is equivalent to
the canonical morphism \ $\O(-1)\otimes D\to\O\otimes\Lambda^2V$. Hence
$\Coker(f)$ depends only on $\ell$, and we can denote $E_\ell=\Coker(f)$. We
have proved that $Y$ is canonically isomorphic to $\P(V^*)$. Using the canonical
complex
\[0\lra\O(-2)\lra\O(-1)\otimes V\lra\O\otimes\Lambda^2V\lra\O(1)\otimes
\Lambda^3V\lra 0\]
it is easy to see that we have an exact sequence
\[0\lra\O_\ell(-1)\lra E_\ell\lra\O(1)\lra 0 .\]
We have $\text{Ext}^1(\O(1),\O_\ell(-1))\simeq\C$ and the preceeding extension
is not trivial, because $E_\ell$ is simple.

A morphism $\O(-2)\to E_\ell$ is non injective if and only if its image is
contained in $\O_\ell(-1)$. It follows that the fiber over $\ell$ of the
projection $\widetilde{Y}\to Y$ is precisely \ $\P(\H^0(E_\ell(2))\backslash
\H^0(\O_\ell(1))$, and that \ $(\boldsymbol W/G)\backslash X_0(4,3)$ is
canonically isomorphic to the projective bundle $S^2(T_{\P(V^*)}(1))$ over
$\P(V^*)$.

The closed subvariety $\widetilde{Y}$ can also be described precisely. It
consists of the non trivial extensions
\[0\lra\O_\ell(-1)\lra\F\lra\O_X(1)\lra 0\]
where $\ell$ is a line and $X$ a cubic. We have
$\dim(\text{Ext}^1(\O_X(1),\O_\ell(-1))=3$, and $\widetilde{Y}$ is a projective
bundle over \ $\P(V^*)\times\P(S^3V^*)$ . The generic sheaves in
$\widetilde{Y}$ are obtained as follows: take $\ell$ and $X$ transverse, the
sheaves $\F$ are obtained by glueing $\O_\ell(2)$ and $\O_X(1)$ at the
intersection points of $\ell$ and $X$.
\end{subsub}

\SEP

We will now describe the closed stratum $X_1(4,1)$. Using the description of
$\mathbb{W'}$ in \ref{clos1} we get easily

\SEP

\begin{subsub}\label{3.A4}{\bf Proposition : }The sheaves of $X_1(4,1)$ are the
kernels of the
surjective morphisms $\O_C(1)\to\O_P$, $C$ beeing a quartic curve
in $\P^2$ and $P$ a closed point of $C$.
\end{subsub}

\end{sub}

\sepsec

\section{Euler Characteristic Two}

\Ssect{Preliminaries}{prel4}

We quote 4.5 from \cite{maican}:

\SEP

\begin{subsub}\label{4.A}{\bf Theorem : }
Let $\F$ be a sheaf in $\MM(4,2)$ satisfying $h^0(\F(-1))=0$ and $h^1(\F)=0$.
Then $h^0(\F\tensor \Om^1(1))$ is zero or one. The sheaves of the first kind are
precisely the sheaves with resolution of the form
\bdm
\tag{i}
0 \lra 2\O(-2) \lra 2\O \lra \F \lra 0.
\edm
The sheaves of the second kind are precisely the sheaves with resolution
\bdm
\tag{ii}
0 \lra 2\O(-2) \oplus \O(-1) \stackrel{\f}{\lra} \O(-1) \oplus 2\O \lra \F \lra 
0,
\edm
\bdm
\f= \begin{pmatrix}
X_1 & X_2 & 0 \\
\star & \star & Y_1 \\
\star & \star & Y_2
\end{pmatrix}
.
\edm
Here $X_1,X_2 \in V^*$ are linearly independent one-forms and the same for
$Y_1,Y_2 \in V^*$.
\end{subsub}

\SEP

\begin{subsub}\label{4.A1}{\bf Remark : } \rm Let $\E$ be a coherent sheaf on
$\P^2$ with Hilbert polynomial $4n+2$. Then there $\E$ is isomorphic to the
cokernel on an injective morphism
\[f:2\O(-2)\lra 2\O\]
if and only if \ $h^0(\F\tensor \Om^1(1))=0$. We have then also $h^0(\F(-1))
=h^1(\F)=0$. This can be seen with a {\em generalized Beilinson spectral
sequence} using exceptional bundles on $\P^2$ (cf. \cite{dr1},
\cite{go_ru},\cite{dr3} 5-). Here we use the triad $(\O(-2),\O,\Omega^1(2))$
(cf. \cite{dr1}) instead of $(\O(-2),\O(-1),\O)$. Then $f$ is equivalent to the
canonical morphism
\[\O(-2)\otimes\H^1(\E(-1))\lra\O\otimes\H^1(\E\otimes Q_2(-1))\]
(where $Q_2$ is the exceptional bundle cokernel of the canonical morphism \
\hfil\break$\O(-2)\to S^2V\otimes\O$).
\end{subsub}

\SEP

\begin{subsub}\label{4.B}{\bf Theorem : }
There are no sheaves $\F$ in $\MM(4,2)$ satisfying the relations
\hfil\break$h^0(\F(-1))=0$ and $h^1(\F)>0$.
\end{subsub}

\Proof Let $\F$ be a semi-stable sheaf with Hilbert polynomial $4m+2$.
Let \hfil\break\mM{p=h^0(\F\ot\Omega^1(1))} and suppose that \ $h^1(\F)>0$ \ and
\ $h^0(\F(-1))=0$. Then we can write \ $h^1(\F)=q+1$ , with $q\geq 0$. The
Beilinson diagram (2.2.3) of $\F$ is
\[\xymatrix{\O(-2)\ot\C^2 & & \O(-1)\ot\C^p & & \O\ot\C^{q+1} \\
0 & & \O(-1)\ot\C^p & & \O\ot\C^{q+3}
}\]
The morphism \ $\varphi_2:\O(-1)\ot\C^p\to\O\ot\C^{q+1}$ \ is surjective,
hence we must have $p\geq q+3$.

The morphism \ $\varphi_3:\O(-1)\ot\C^p\to\O\ot\C^{q+3}$ \ is injective,
hence we have \ $p\leq q+3$. Finally we get $p=q+3$.

The Beilinson complex of $\F$ is then
\[0\lra
2\O(-2)\oplus(q+3)\O(-1)\lra(q+3)\O(-1)\oplus(q+3)\O\lra(q+1)\O\lra 0 .\]
Let $\E=\Coker(\varphi_3)$ and $E=\ker(\varphi_2)$ (it is a rank 2 vector
bundle). We have a commutative diagram with exact rows
\[\xymatrix{0\ar[r] & (q+3)\O(-1)\ar[r]^-{\varphi_3}\flinc[d] &
(q+3)\O\ar[r]\flinc[d] & \E\ar[r]\ar[d]^\alpha & 0\\
0\ar[r] & 2\O(-2)\oplus(q+3)\O(-1)\ar[r] & E\oplus(q+3)\O\ar[r] & \F\ar[r] &
0 }\]
If follows that $\ker(\alpha)\subset\ 2\O(-2)$. But $\E$ is a torsion sheaf.
Hence $\alpha$ is injective.
But the Hilbert polynomial of $\E$
is $(q+3)(m+1)$. This contradicts the semi-stability of $\F$.
\carre

\end{sub}

\sepsub

\Ssect{The two sub-strata of the open stratum}{open2}

Let $X$ be the open subset of $\M(4,2)$ corresponding to sheaves $\F$ satisfying
the conditions
\bdm
h^0(\F(-1))=0, \qquad h^1(\F)=0, \qquad  h^0(\F \tensor\Om^1(1))\leq 1 .
\edm
It is the disjoint union of the open subset $X_0$ and the locally closed subset
$X_1$, where for $i=0,1$, $X_i$ corresponds to sheaves $\F$ such that \
$ h^0(\F \tensor\Om^1(1))=i$ .
We will first describe the subsets $X_0$, $X_1$.

\SEP

\begin{subsub}\label{4.D0} {\bf The open subset $X_0$ -} \rm Let
$W=\Hom(2\O(-2),2\O)$ on which acts the reductive group
$G_0=\text{GL}(2)\times\text{GL}(2)$. Let $W_0\subset W$ be the set of
morphisms from \ref{4.A}(i) that is, the set of injective morphisms
\bdm
\f : 2\O(-2) \lra 2\O.
\edm
The corresponding $S^2V$-Kronecker modules \ $S^2V\otimes\C^2\to\C^2$ \ are
then semistable (cf. \ref{KM}). Hence $W_0$ is an open $G_0$-invariant subset
inside the set $W^{ss}$ of semistable points and contains the set
$W^{ss}\setminus W^s$ of properly semistable points. $W^{ss}$ is the set of $2
\times 2$-matrices with entries in $S^2V^*$ having linearly independent rows and
columns. $W^{ss} \setminus W^s$ is the subset of matrices equivalent, modulo row
and column operations, to a matrix having a zero entry. Incidentally, note that
$W_0$ is a proper subset of $W^{ss}$ (cf. \ref{desc2}). Thus $W_0/\!\!/G_0$ is a
proper open subset of the projective variety $N(6,2,2)$.

The morphism \ $\rho_0: W_0 \to X_0$ \ given by \ $\xymatrix{\f \fmaps[r] &
[\Coker(\f)]}$ \ is $G_0$-invariant. We claim that $\rho_0(\f_1)=\rho_0(\f_2)$
if and only if \ $\overline{G_0\f_1} \cap \overline{G_0\f_2}\neq \emptyset$.
This is clear if $\F=\Coker(\f_1)$ is stable, in fact, in this case, we have
$G_0\f_1 = G_0\f_2$. If $\F$ is properly semistable, then there is an extension
\bdm
0 \lra \O_{C_1} \lra \F \lra \O_{C_2} \lra 0
\edm
for some conics $C_1= \{ f_1=0 \}$ and $C_2= \{ f_2 = 0 \}$.
From the horseshoe lemma we get a resolution
\bdm
0 \lra 2\O(-2) \stackrel{\psi}{\lra} 2\O \lra \F \lra 0, \qquad \psi =
\begin{pmatrix}
f_2 & 0 \\
\star & f_1
\end{pmatrix}
.
\edm
Thus $f_1 \oplus f_2$ is in the closure of the orbit $G_0\psi = G_0\f_1$.
Analogously, $f_1 \oplus f_2$ is in the closure of $G_0\f_2$.
\end{subsub}

\SEP

\begin{subsub}\label{4.D}{\bf Theorem : }
The good quotient $W_0/\!\!/G_0$ is isomorphic to $X_0$.
The subvariety \mM{(W^{ss}\setminus W^s)/\!\!/G_0} \ of
$W_0/\!\!/G_0$ given by properly semistable points is isomorphic to the
subvariety of $\MM(4,2)$ given by properly semistable sheaves and is isomorphic
to the symmetric space \ $\big(\P(S^2V^*) \times \P(S^2
V^*)\big)/\text{\emph{S}}_2$.
\end{subsub}

\Proof We will show that $\rho_0$ is a categorical quotient map and the
isomorphism $W_0/\!\!/G_0\isom X_0$ will follow from the uniqueness of the
categorical quotient. Let \hbox{$f: W_0 \to Y$} \ be a $G_0$-invariant morphism
of varieties. On the closure of each $G_0$-orbit $f$ is constant hence, by the
above remark, $f$ is constant on the fibers of $\rho$. Thus $f$ factors through
a map \ $g:X_0 \to Y$.

We continue the proof as at \ref{3.B}. We need to construct resolution
\ref{4.A}(i)
starting from the Beilinson spectral sequence for $\F$. Tableau (2.2.3) takes
the form
\bdm
\xymatrix
{
2\O(-2) & 0 & 0 \\
0 & 0 & 2\O
}
\edm
and (2.2.5) yields the resolution
\bdm
0 \lra 2\O(-2) \stackrel{\f_5}{\lra} 2\O \lra \F \lra 0.
\edm
This allows us to construct the morphism \ $\xi : E \to W_0$ as at \ref{3.B} \
and proves the isomorphism $W_0/\!\!/G_0\isom X_0$.

The statement about properly semistable sheaves follows from the easily
verifiable fact that $\rho_0(\f)$ is properly semistable if and only if $\f$ is
properly semistable. \carre

\SEP

\begin{subsub}\label{4.E0} {\bf The locally closed subset $X_1$ -} \rm
Let $W_1$ be the set of morphisms from \ref{4.A}(ii) that is, the set of
injective morphisms
\bdm
\f : 2\O(-2) \oplus \O(-1) \lra \O(-1) \oplus 2\O
\edm
satisfying $\f_{12}=0$ and such that $\Coker(\f)$ is semistable.
The linear algebraic group
$$G=\Aut(2\O(-2) \oplus \O(-1)) \times \Aut(\O(-1) \oplus 2\O)$$
acts on $W_1$ in an obvious way. Let $X_1$ be the locally closed subset of
$\M(4,2)$ defined by the conditions
\bdm
h^0(\F(-1))=0, \qquad h^1(\F)=0, \qquad h^0(\F \tensor \Om^1(1))=1.
\edm
We equip $X_1$ with the canonical induced reduced structure. The morphism
$$\xymatrix@R=4pt{\rho_1: W_1 \ar[r] & X_1\\ \f \fmaps[r] & [\Coker(\f)]}$$
is surjective and its fibers are $G$-orbits.
\end{subsub}

\SEP

\begin{subsub}\label{4.E}{\bf Theorem : }
The morphism $\rho_1:W_1\to X_1$ is a geometric quotient by $G$.
\end{subsub}

\Proof According to \cite{mumf}, prop. 0.2, since the fibers of $\rho_1$ are the
$G$-orbits we need only to prove that $\rho_1$ is a categorical quotient.
As in the proof of \ref{3.B}, we need to recover resolution
\ref{4.A}(ii) of a sheaf $\F$ in $X_1$ starting from the Beilinson spectral
sequence for $\F$. Tableau (2.2.3) is
\bdm
\xymatrix
{
2\O(-2) \ar[r]^{\f_1} & \O(-1) & 0 \\
0 & \O(-1) \ar[r]^{\f_4} & 2\O
}.
\edm
$\Coker(\f_1)$ cannot be of the form $\O_L$ for a line $L \subset \P^2$ because,
by semistability, $\F$ cannot surject onto such a sheaf. Thus $\Coker(\f_1)$ is 
supported on a point and $\Ker(\f_1) \isom \O(-3)$.
Clearly $\f_5$ lifts to a morphism $\psi_5: \O(-3) \lra 2\O$. We have
a resolution
$$\xymatrix{0\ar[r] & \O(-3) \oplus \O(-1)\ar[rr]^-{[\psi_5, \f_4 ]} & &
2\O\ar[r] & \Coker(\f_5)\ar[r] & 0 \ .}$$
We now apply the horseshoe lemma to the extension
\bdm
0 \lra \Coker(\f_5) \lra \F \lra \Coker(\f_1) \lra 0,
\edm
to the resolution of $\Coker(\f_5)$ given above and to the resolution
\bdm
0 \lra \O(-3) \lra 2\O(-2) \lra \O(-1) \lra \Coker(\f_1) \lra 0.
\edm
We obtain the exact sequence
\bdm
0 \lra \O(-3) \lra \O(-3) \oplus \O(-1) \oplus 2\O(-2) \lra 2\O \oplus \O(-1) 
\lra \F \lra 0.
\edm
From the fact that $h^1(\F)=0$, we see that $\O(-3)$ can be cancelled to yield a
resolution
as in \ref{4.A}(ii). \carre
\end{sub}

\sepsub

\Ssect{The open stratum}{open3}

Let
\[\W = \Hom(2\O(-2)\oplus\O(-1),\O(-1)\oplus 2\O) .\]
The elements of $\W$ are represented as matrices
\[\begin{pmatrix}X_1 & X_2 & \alpha\\ q_{11} & q_{12} & Y_1\\
q_{21} & q_{22} & Y_2\end{pmatrix} ,\]
where $\alpha\in\C$, $X_1$, $X_2$, $Y_1$, $Y_2$ are linear forms and
$q_{11}$,$q_{12}$,$q_{21}$,$q_{22}$ quadratic forms on $V$. The linear algebraic
group
$$G=\Aut(2\O(-2) \oplus \O(-1)) \times \Aut(\O(-1) \oplus 2\O)$$
acts on $\W$ as in \ref{4.E0} and $W_1$ is a locally closed $G$-invariant
subset of $\W$.

We are in  the situation of \ref{MSM}, and we will use the polarization
\[\sigma=(\frac{1-\mu}{2},\mu,\mu,\frac{1-\mu}{2}) ,\]
where $\mu$ is a rational number such that $\frac{1}{3}<\mu<\frac{1}{2}$.

\SEP

\begin{subsub}\label{4.E1}{\bf Lemma : } let $f\in\W$ be an injective morphism.
Then
$\Coker(f)$ is semi-stable if and only if $f$ is $G$--semi-stable with respect
to $\sigma$.
\end{subsub}

\Proof
 Let $m_1,m_2,n_1,n_2$ be integers such that $0\leq m_1\leq 2$, $0\leq m_2\leq
1$, $0\leq n_1\leq 1$, $0\leq n_2\leq 2$. Let $f\in\W$. We say that \ {\em
``$(m_1,m_2)\to(n_1,n_2)$ is forbidden for $f$ ''} if for every linear subspaces
$M_1\subset\C^2$, $M_2\subset\C$, $N_1\subset\C$, $N_2\subset\C^2$ such that
$\dim(M_i)=m_i$  and $\dim(N_j)=n_j$ we don't have
\[f((\O(-2)\otimes M_1)\oplus(\O(-1)\otimes M_2))\subset(\O(-1)\otimes N_1)
\oplus(\O\otimes N_2) .\]
Then $f$ is $G$--semi-stable with respect to $\sigma$ if and only the following
are forbidden for $f$ :
\[(2,0)\to(1,0), (2,0)\to(0,1), (0,1)\to(0,1), (1,1)\to(0,2), (1,1)\to(1,0),\]
and the cases where $n_1=n_2=0$ and $m_1+m_2\not=0$. The results follows then
easily from \ref{open2}.
\carre

\SEP

The elements of the group $G$ can be seen as pairs of matrices
\[(\nu_1,\nu_2) \ = \ \biggl(\begin{pmatrix}\alpha & 0\\ \phi & A\end{pmatrix}
, \begin{pmatrix}B & 0\\ \psi & \beta\end{pmatrix}\biggr)\]
where $A,B\in\GL(2)$, $\alpha,\beta\in\C^*$, $\phi$ is a column vector
$\begin{pmatrix}\phi_1\\ \phi_2\end{pmatrix}$ and $\psi=(\psi_1,\psi_2)$ a pair
of linear forms on $V$. As in \ref{4.D0} we consider the space \
$W=\Hom(2\O(-2),2\O)$ \ on which acts the reductive group \
$G_0=\GL(2)\times\GL(2)$ . Let
\[\tau:G\lra G_0\]
(cf. \ref{4.D0}) be the morphism of groups defined by \ $\tau(\nu_1,\nu_2)=
\alpha\beta^{-1}AB^{-1}$. An easy calcultation shows that

\SEP

\begin{subsub}\label{4.E2}{\bf Lemma : } The morphism
$$\Delta:\W\lra W$$
defined by
\[\Delta\biggl(\begin{pmatrix}X_1 & X_2 & \alpha\\ q_{11} & q_{12} & Y_1\\
q_{21} & q_{22} & Y_2\end{pmatrix}\biggr) \ = \
\alpha\begin{pmatrix}q_{11} &
q_{12}\\q_{21} & q_{22}\end{pmatrix}-\begin{pmatrix}Y_1\\Y_2\end{pmatrix}
(X_1,X_2)\]
is compatible with $\tau$, i.e. for every $w\in\W$ and $g\in G$ we have
$\Delta(gw)=\tau(g)\Delta(w)$.
\end{subsub}

\SEP

Note that the image of $\Delta$ is $W^{ss}$.

The locus of non injective morphism in N(6,2,2) is isomorphic to \
$\P^2\times\P^2$. To a pair $(P_0,P_2)$ of points of $\P^2$ corresponds the
$\text{GL}(2)\times\text{GL}(2)$-orbit of the matrix
$\begin{pmatrix}\alpha_0\beta_0 & \alpha_0\beta_1\\ \alpha_1\beta_0 &
\alpha_1\beta1\end{pmatrix}$, where for $i=0,1$, $(\alpha_i,\beta_i)$ is a
pair of linear forms defining $P_i$. Thus we have an isomorphism
\[X_0 \ \simeq \ N(6,2,2)\backslash(\P^2\times\P^2) .\]
Let \ $\pi:W^{ss}\to N(6,2,2)$ \ be the quotient morphism. We have
\[\Delta(W_1)=\pi^{-1}(\P^2\times\P^2)\subset W^s .\]

We have a surjective $G$--invariant morphism
$$\xymatrix@R=4pt{\rho: \W \ar[r] & X\\ \f \fmaps[r] & [\Coker(\f)]}$$

\SEP

\begin{subsub}\label{4.E3}{\bf Theorem : }
The morphism $\rho$ is a good quotient by $G$.
\end{subsub}

\Proof Let $\W_0=\rho^{-1}(W^s)$. According to \cite{mumf}, prop. 0.2, the
restriction of $\rho$, $\W_0\to W^s$ is a geometric quotient by $G$. Let \
\[\W_2=\Delta^{-1}\big(N(6,2,2)\backslash(\P^2\times\P^2)\big)=\rho^{-1}(X_0)
.\]
We will show that the restriction of $\rho$, $\W_2\to X_0$ is a good quotient.
Since \ \mM{\W=\W_0\cup\W_2} \ it follows easily from definition \ref{2.L0}
that $\rho$ is a good quotient.

Let \ $H=\C^*\times(V^*)^4$.
We have an isomorphism
\[\theta:\pi^{-1}(N(6,2,2)\backslash(\P^2\times\P^2)\big)\times H\to\W_2\]
given by
\[\theta(q_0,\mu,X_1,X_2,Y_1,Y_2)=\begin{pmatrix}\vec X & \mu\\Q & \vec Y
\end{pmatrix} \ ,
\]
with
\[{\vec X}=(X_1,X_2) , \quad {\vec Y}=\begin{pmatrix}Y_1\\Y_2\end{pmatrix} ,
\quad Q=\frac{1}{\mu}\big(q_0+\begin{pmatrix}Y_1\\Y_2\end{pmatrix}(X_1,X_2)
\big) .\]
It is easy to verify that the $G$--orbits of $\W_2$ are of the form $Y\times H$,
where $Y$ is a $G_0$--orbit of $N(6,2,2)\backslash(\P^2\times\P^2)$. The fact
that $\rho$ is a good quotient by $G$ follows then immediately from the fact
that the restriction of $\pi$,
\[\pi^{-1}(N(6,2,2)\backslash(\P^2\times\P^2)\big)\lra
N(6,2,2)\backslash(\P^2\times\P^2)\]
is a good quotient by $G_0$.
\carre

\SEP

From $\Delta$ and \ref{4.E3} we get a surjective morphism
\[\delta:X\lra N(6,2,2)\]
which induces an isomorphism \ $X_0\simeq N(6,2,2)\backslash(\P^2\times\P^2)$.

\SEP

\begin{subsub}\label{4.E4}\bf The fibers of $\delta$ over points of
$\P^2\times\P^2$ - \rm Let $P_1$, $P_2$ distinct points of $\P^2$. Let $X_i$,
$Z$ ($i=1,2$) be linearly independant linear forms on $V$ vanishing at $P_i$.
Then $\delta^{-1}(P_1,P_2)$ contains the sheaves which are cokernels of
injective morphisms of type
\[\begin{pmatrix}X_1 & Z & 0\\ q_{11} & q_{12} & Z\\ q_{21} & q_{22} &
X_2 \end{pmatrix} .\]
Each $G$--orbit of such a morphism contains a matrix of the following type :
\[\begin{pmatrix}X_1 & Z & 0\\ \alpha X_2^2 & q_{12} & Z\\ q_{21} & \beta X_1^2
& X_2 \end{pmatrix}\]
where the quadratic form $q_{12}$ in $X_1$, $X_2$, $Z$ has no term in $Z^2$,
and all the matrices of this type in the $G$--orbit are obtained by replacing
the submatrix $\begin{pmatrix}\alpha X_2^2 & q_{12}\\q_{21} & \beta X_1^2
\end{pmatrix}$ by a nonzero multiple. The corresponding morphism is non
injective if and only if the submatrix vanishes. It follows that \
$\delta^{-1}(P_1,P_2)\simeq\P^{12}$ .

Let $P_1$ be a point of $\P^2$, and $X_1$, $Z$ be linearly independant linear
forms on $V$ vanishing at $P_1$. Let $X_2$ be a linear form such that
$(X_1,X_2,Z)$ is a basis of $V^*$. Then $\delta^{-1}(P_1,P_1)$ contains the
sheaves which are cokernels of injective morphisms of type
\[\begin{pmatrix}X_1 & Z & 0\\ q_{11} & q_{12} & Z\\ q_{21} & q_{22} &
X_1 \end{pmatrix} .\]
Each $G$--orbit of such a morphism contains a matrix of the following type :
\[\begin{pmatrix}X_1 & Z & 0\\ \alpha X_2^2 & q_{12} & Z\\ q_{21} & \beta X_2^2
& X_1 \end{pmatrix}\]
where the quadratic form $q_{12}$ in $X_1$, $X_2$, $Z$ has no term in $Z^2$,
and all the matrices of this type in the $G$--orbit are obtained by replacing
the submatrix $\begin{pmatrix}\alpha X_2^2 & q_{12}\\q_{21} & \beta X_2^2
\end{pmatrix}$ by a nonzero multiple. The corresponding morphism is non
injective if and only if the submatrix vanishes or is a multiple of
$\begin{pmatrix}X_2^2 & 0\\0 & -X_2^2\end{pmatrix}$. Hence we see that
$\delta^{-1}(P_1,P_1)$ is isomorphic to the complement of a point in $\P^{14}$.
\end{subsub}

\SEP

\begin{subsub}\label{4.E5}{\bf Theorem : } Let $\widetilde{\mathbf{N}}$ be the
blowing-up of
$N(6,2,2)$ along $\P^2\times\P^2$. Then $X$ is isomorphic to an open subset of
$\widetilde{\mathbf{N}}$.
\end{subsub}

\Proof We have \ $\delta^{-1}(\P^2\times\P^2)=W_1$ , which is a smooth
hypersurface of $X$. It follows from the universal property of the blowing-up
that $\delta$ factors through $\widetilde{\mathbf{N}}$ : we have a morphism \
$\widetilde{\delta}:X\to\widetilde{\mathbf{N}}$ \ such that \
$\delta=p\circ\widetilde{\delta}$, $p$ beeing the projection \
$\widetilde{\mathbf{N}}\to N(6,2,2)$. We want to prove that
$\widetilde{\delta}$ induces an isomorphism $X\simeq\widetilde{\delta}(X)$. For
this it suffices to prove that $\widetilde{\delta}$ does not contract $W_1$ to
a subvariety of codimension $\geq 2$ (cf. \cite{sha}, II,4, theorem 2).

The morphism $\widetilde{\delta}$ can be described precisely on $W_1$. Recall
that $p^{-1}(\P^2\times\P^2)$ is the projective bundle
$\P(\N_{\P^2\times\P^2})$, $\N_{\P^2\times\P^2}$ beeing the normal bundle of
$\P^2\times\P^2$ in $N(6,2,2)$. Let $w\in W_1$ and $C$ a smooth curve through
$w$ in $X$, not tangent to $W_1$ at $w$. Then $\widetilde{\delta}(w)$ is the
image in $\N_{\P^2\times\P^2,\delta(w)}$ of the tangent to $C$ at $w$. Suppose
first that $\delta(w)$ is a pair of distinct points $(P_1,P_2)$ and that
$w$ is the cokernel of a morphism
\[\phi_0=\begin{pmatrix}X_1 & Z & 0\\ q_{11} & q_{12} & Z\\ q_{21} & q_{22} &
X_2 \end{pmatrix}\]
(we use the notations of \ref{4.E4}). Then \
$\overline{\phi_0}=\begin{pmatrix}ZX_1 & Z^2\\ X_1X_2 & ZX_2\end{pmatrix}$ \ is
a point of $W^s$ over $w$ and $\N_{\P^2\times\P^2,\delta(w)}$ is
isomorphic to $\N_{\Gamma,\overline{\phi_0}}$, where $\Gamma\subset W^s$ is the
inverse image of $\P^2\times\P^2$ and $\N_\Gamma$ its normal bundle. The vector
space $\N_{\Gamma,\overline{\phi_0}}$ is a quotient of $W$. Supppose that $C$ is
defined by the family $(\phi_t)$ (for $t$ in a neighbourhood of $0$ in $\C$),
with
\[\phi_t=\begin{pmatrix}X_1 & Z & t\\ q_{11} & q_{12} & Z\\ q_{21} & q_{22} &
X_2 \end{pmatrix} \ .\]
Then from the formula defining $\Delta$ we deduce that $\widetilde{\delta}(w)$
is the image of $\begin{pmatrix}q_{11} & q_{12}\\ q_{21} & q_{22}\end{pmatrix}$
in $\N_{\Gamma,\overline{\phi_0}}$. It follows from \ref{4.E4} that
$\N_{\P^2\times\P^2,\delta(w)}$ is contained in the image of
$\widetilde{\delta}$, and that $\widetilde{\delta}$ does not contract $W_1$ to
a subvariety of codimension $\geq 2$.
\carre

\SEP

\begin{subsub}\label{4.E6}{\bf Remark : }\rm Using \ref{4.E4} and the
proof of \ref{4.E5} is is not difficult to prove that \
$\widetilde{\mathbf{N}}\backslash X\simeq\P^2$ .
\end{subsub}

\end{sub}

\sepsub

\Ssect{The closed stratum}{clos2}

\begin{subsub}\label{4.F}{\bf Theorem : }
A sheaf $\F$ giving a point in
$\MM(4,2)$ satisfies $h^0(\F(-1))> 0$
if and only if there is a quartic $C \subset \P^2$ such that $\F \isom \O_C(1)$.
The subvariety of such sheaves in $\MM(4,2)$ is isomorphic to $\P(S^4V^*)$.
\end{subsub}

\Proof As explained in the comments before \ref{2.C}, a nonzero morphism
\hbox{$\O \to \F(-1)$} \ must factor through an injective morphism \
$\O_C \to \F(-1)$ \ for a curve $C\subset \P^2$. From the semistability of $\F$
we see that $\F$ must be a quartic, and the above morphism must be an
isomorphism. \carre
\end{sub}

\sepsub

\Ssect{Description of the strata}{desc2}

We will describe $X_1$.

Let \ $\f:2\O(-2)\oplus\O(-1)\to\O(-1)\oplus 2\O$ \ be an injective morphism as
in \ref{4.A} (ii) and $\F=\Coker(\f)$. Let $P$ (resp. $Q$) be the point of
$\P^2$ defined by the linear forms $X_1$, $X_2$ (resp. $Y_1$, $Y_2$). We have a
commutative diagram with exact rows
\[\xymatrix{
0\ar[r] & \O(-1)\ar[r]\flinc[d] & 2\O\ar[r]\flinc[d] &
\I_Q(1)\ar[r]\ar[d]^\alpha & 0\\
0\ar[r] & 2\O(-2)\oplus\O(-1)\ar[r]^-\f & \O(-1)\oplus 2\O\ar[r] &
\F\ar[r] & 0
}\]
It follows that we have $\Ker(\alpha)\simeq\O(-3)$ and
$\Coker(\alpha)\simeq\O_P$ (the structural sheaf of $P$). So we have an exact
sequence
\[0\lra\O(-3)\stackrel{\beta}{\lra}\I_Q(1)\lra\F\lra\O_P\lra 0 .\]
Up to a scalar multiple there is only one surjective morphism \
$\O_C(1)\to\O_Q$ . We denote by $\O_C(1)(-Q)$ its kernel. If $C$ is not smooth
at $P$ this sheaf can be non locally free. We have \
$\Coker(\beta)\simeq\O_C(1)(-Q)$, and an exact sequence
\setcounter{subsub}{1}
\begin{equation}\label{4.G2}
0\lra\O_C(1)(-Q)\lra\F\lra\O_P\lra 0 .
\end{equation}
Hence we could denote $\F$ by ``$\O_C(1)(P-Q)$''. This notation is justified if
$P\not=Q$ and if $C$ is smooth at $P$. So we can already state

\SEP

\begin{subsub}\label{4.G}{\bf Proposition : }The generic sheaf in $X_1$ is of
the form
$\O_C(1)(P-Q)$, where $C$ is a smooth quartic and $P$, $Q$ distinct points of
$C$.
\end{subsub}

\SEP

We want to prove now that the notation ``$\O_C(1)(P-Q)$'' has a meaning if
$P\not=Q$ even if $C$ is not necessarily smooth. We have to compute
$\text{Ext}^1_{\O_C}\big(\O_P,\O_C(1)(-Q)\big)$. First we note that
\[\text{Ext}^1_{\O_C}\big(\O_P,\O_C(1)(-Q)\big) \ \simeq \
\text{Ext}^1_{\O_{\P^2}}\big(\O_P,\O_C(1)(-Q)\big) \]
(this can be seen using prop. 2.3.1 of \cite{dr4}). We use now the exact
sequence
\[0\to\O(-3)\to\I_Q(1)\to\O_C(1)(-Q)\to 0\]
and get the exact sequence
\[0\lra\text{Ext}^1_{\O_{\P^2}}\big(\O_Q,\I_P(1)\big)\lra\text{Ext}^1_{\O_{\P^2}
} \big(\O_P,\O_C(1)(-Q)\big)\lra\text{Ext}^2_{\O_{\P^2}}\big(\O_Q,\O(-3)\big)\]
\[\lra\text{Ext}^2_{\O_{\P^2}}(\O_Q,\I_P(1))\lra\text{Ext}^2_{\O_{\P^2}}(\O_P,
\O_C(1)(-Q))\lra 0 .\]
Using the exact sequence \ $0\to\O(-1)\to 2\O\to\I_P(1)\to 0$ , we find that
\[\text{Ext}^1_{\O_{\P^2}}\big(\O_Q,\I_P(1)\big)=\lbrace 0\rbrace \ ,\quad
\text{Ext}^2_{\O_{\P^2}}\big(\O_Q,\I_P(1)\big)\simeq\C\]
if $P\not=Q$ and
\[\text{Ext}^1_{\O_{\P^2}}\big(\O_P,\I_P(1)\big)\simeq\C \ ,\quad
\text{Ext}^2_{\O_{\P^2}}\big(\O_P,\I_P(1)\big)\simeq\C^2 .\]
It now follows easily that \
$\text{Ext}^1_{\O_C}\big(\O_P,\O_C(1)(-Q)\big)\simeq\C$
if \ $P\not=Q$, hence there is only one non trivial extension ({\ref{4.G2}}) and
the notation $\O_C(1)(P-Q)$ is justified in this case.

If $P=Q$ we first remark that $C$ is never smooth at $P$, and this implies
that the morphism
$$\text{Ext}^2_{\O_{\P^2}}\big(\O_Q,\O(-3)\big)\lra\text{Ext}^2_{\O_{\P^2}}
\big(\O_Q ,\I_P(1)\big)$$
vanishes: by Serre duality it is the transpose of
\[\Hom(\I_P(1),\O_Q)\lra\Hom(\O(-3),\O_Q)\]
which is just the multiplication by an equation of $C$. It follows that we have
\hfil\break $\text{Ext}^1_{\O_C}\big(\O_P,\O_C(1)(-P)\big)\simeq\C^2$ \ and an
injective
map
\[\xymatrix{\lambda:\C=\text{Ext}^1_{\O_{\P^2}}\big(\O_P,\I_P(1)\big)\flinc[r] &
\text{Ext}^1_{\O_C}\big(\O_P,\O_C(1)(-P)\big)\simeq\C^2 .}\]
Its image corresponds to the extension ({\ref{4.G2}}) given by $\F=\O_C(1)$ and
the other extensions, which are in $X_1$ are defined by the other elements of
$\text{Ext}^1_{\O_C}\big(\O_P,\O_C(1)(-P)\big)$.
\end{sub}

\sepsec

\section{Euler Characteristic Four}

\Ssect{The open stratum}{open4}

Let $X_0$ be the open subset of $\M(4,4)$ corresponding to sheaves $\F$ such
that \hfil\break$h^0(\F(-1))=0$. The complement $\M(4,4)\backslash X_0$ is the
{\em theta divisor} (cf. \cite{lepotier}). Combining \cite{lepotier}, 4.3 and
\cite{dr2}, th\'eor\`eme 2, we obtain

\SEP

\begin{subsub}{\bf Theorem : }
1 - The sheaves on $\P^2$ with Hilbert polynomial $4m+4$ satisfying
$h^0(\F(-1))=0$ are precisely the sheaves which are isomorphic to the cokernel
of an injective morphism
\[f:4\O(-1)\lra 4\O \ .\]
Moreover, $\F$ is not stable if and only if $\f$ is equivalent,
modulo operations on rows and columns, to a morphism of the form
\bdm
\begin{pmatrix}
\f_{11} & 0 \\
\f_{21} & \f_{22}
\end{pmatrix}
\qquad \text{\emph{with}} \quad \f_{22} : m\O(-1) \lra m\O, \quad m=1,
\, 2 \ \text{\emph{or}}\ 3.
\edm

2 - Let $\mathbf{N}_0$ denote the open subset of $N(3,4,4)$ corresponding to
injective morphisms. By associating $\F$ to the
$\big(\GL(4)\times\GL(4)\big)$-orbit of $\f$ we get an isomorphism \
$\mathbf{N}_0\simeq X_0$.
\end{subsub}

\SEP

According to \cite{lepotier} the complement $N(3,4,4)\backslash X_0$ is
isomorphic to $\P^2$, the inclusion \hfil\break$X_0\subset N(3,4,4)$ can be
extended to a morphism \ $\M(4,4)\to N(3,4,4)$ \ which is the blowing-up of
$N(3,4,4)$ along $\P^2$.

\end{sub}

\sepsub

\Ssect{The closed stratum}{clos4}

\begin{subsub}\label{5.C}{\bf Proposition : }
The sheaves $\F$ in $\MM(4,4)$ satisfying $h^0(\F(-1))=1$ are
precisely the sheaves with resolution
\bdm
0 \lra \O(-2) \oplus \O(-1) \stackrel{\f}{\lra} \O \oplus \O(1) \lra \F \lra 0, 
\qquad \f_{12}\neq 0.
\edm
Moreover, $\F$ is not stable if and only if $\f_{12}$ divides
$\f_{11}$ or $\f_{22}$.
\end{subsub}

\Proof We assume that $\F$ has a resolution as above and we need to
show that $\F$ is semistable. Assume that there is a destabilizing subsheaf $\E$
of $\F$. We may assume that $\E$ is semistable. As $\F$ is generated by global
sections, we must have $h^0(\E) < h^0(\F)=4$. Thus $\E$ is in $\M(2,3)$,
$\M(1,3)$ or in $\M(1,2)$. In other words, $\E$ must be isomorphic to $\O_C(1)$,
$\O_L(2)$ or $\O_L(1)$ for a conic $C$ or a line $L$ in $\P^2$. In the first
case we have a commutative diagram
\bdm
\xymatrix
{
0 \ar[r] & \O(-1) \ar[r] \ar[d]^{\b} & \O(1) \ar[r] \ar[d]^{\a} & \E \ar[r] 
\ar[d] & 0 \\
0 \ar[r] & \O(-2) \oplus \O(-1) \ar[r] & \O \oplus \O(1) \ar[r] &\F \ar[r] & 0
}.
\edm
Notice that $\a$ is injective, because it is injective on global sections, hence
$\b$ is injective too. We obtain $\f_{12}=0$, which contradicts the hypothesis
on $\f$. The other two cases similarly lead to contradictions: if $\F \isom
\O_L(2)$, then $\a=0$; if $\F \isom \O_L(1)$, then $\b=0$, which is impossible.

Conversely, we are given $\F$ in $\M(4,4)$ satisfying the condition
$h^0(\F(-1))=1$ and we need to construct a resolution as above. As explained in
the comments above \ref{2.C}, there is an injective morphism \ $\O_C \to \F(-1)$
\ for a curve $C$ in $\P^2$. From the semistability of $\F(-1)$ we see that $C$
is a cubic or a quartic. Assume that $C$ is a cubic curve. The quotient
$\F/\O_C(1)$ has Hilbert polynomial $P(t)=t+1$ and has no zero-dimensional
torsion. Indeed, if $\F/\O_C(1)$ had zero-dimensional torsion $\T \neq 0$, then
the preimage of $\T$ in $\F$ would be a subsheaf which violates the
semistability of $\F$. We conclude that $\F/\O_C(1) \isom \O_L$ for a line $L
\subset \P^2$. We apply the horseshoe lemma to the extension
\bdm
0 \lra \O_C(1) \lra \F \lra \O_L \lra 0
\edm
and to the standard resolutions for $\O_C(1)$ and $\O_L$ to get the desired 
resolution for $\F$. This also proves the statement about properly semistable
sheaves from the claim. Indeed, if $\F$ is properly semistable, then $\F$ has
stable factors $\O_C(1)$ and $\O_L$ and we can apply the horseshoe lemma as
above. For the rest of the proof we may assume that $\F$ is stable.

Next we examine the situation when $C$ is a quartic. Notice that $\F/\O_C(1)$ 
has zero-dimensional support and Euler characteristic 2. There is a subsheaf
$\T \subset \F/\O_C(1)$ satisfying $h^0(\T)=1$. Let $\E$ be the
preimage of $\T$ in $\F$ and $\G=\F/\E$. Notice that $\E$ is semistable because
any subsheaf ruining the semistability of $\E$ must contradict the stability of
$\F$. From the semicontinuity theorem III 12.8 in \cite{hartshorne} we have
$h^0(\F\tensor \Om^1(1)) \ge 4$. Indeed, this inequality
is true for all sheaves in the open dense subset $X_0$ of $\M(4,4)$, so it must
be true for all sheaves giving a point in $\M(4,4)$. We have
\bdm
4 \le h^0(\F \tensor \Om^1(1)) \le h^0(\E \tensor \Om^1(1)) + h^0(\G \tensor 
\Om^1(1)) = h^0(\E \tensor \Om^1(1)) + 1.
\edm
From our results on in section 3, and from the duality theorem \ref{2.A},
we see that there is a resolution
\bdm
0 \lra 2\O(-2) \lra \O(-1) \oplus \O(1) \lra \E \lra 0.
\edm
We apply the horseshoe lemma to the extension
\bdm
0 \lra \E \lra \F \lra \G \lra 0,
\edm
to the resolution of $\E$ from above and to the resolution
\bdm
0 \lra \O(-2) \lra 2\O(-1) \lra \O \lra \G \lra 0.
\edm
The morphism \ $\O \lra \G$ \ lifts to a morphism \ $\O \lra \F$ \ because
$h^1(\E)=0$. We obtain a resolution
\bdm
0 \lra \O(-2) \lra 2\O(-1) \oplus 2\O(-2) \lra \O \oplus \O(-1) \oplus \O(1) 
\lra \F \lra 0.
\edm
From the condition $h^0(\F(-1))=1$ we see that the map $\O(-2) \lra 2\O(-2)$
from the above sequence is nonzero. We may cancel $\O(-2)$ to get the exact 
sequence
\bdm
0 \lra \O(-2) \oplus 2\O(-1) \lra \O(-1) \oplus \O \oplus \O(1) \lra \F \lra 0.
\edm
The morphism \ $2\O(-1)\to\O(-1)$ \ from the above sequence is nonzero,
otherwise $\F$ would surject onto a sheaf of the form $\O_L(-1)$ for a line $L
\subset \P^2$, in violation of the semistability of $\F$. Thus we may cancel
$\O(-1)$ in the above sequence to get the resolution of $\F$ from the claim.
\carre

\SEP

Let $W_1$ be the set of morphisms from \ref{5.C} that is, the set of injective
morphisms
\bdm
\f: \O(-2) \oplus \O(-1) \lra \O \oplus \O(1),
\edm
for which $\f_{12} \neq 0$.
The linear nonreductive algebraic group
$$G=\big(\Aut(\O(-2) \oplus \O(-1)) \times\Aut(\O \oplus \O(1))\big)/\C^*$$
acts on $W_1$ by conjugation. According to \cite{drezet-trautmann}, there is a
good and projective quotient $W_1/\!\!/G$ which contains a geometric quotient as
a proper open subset.For a polarization $\sigma=(\l_1,\l_2,\m_1,\m_2)$
satisfying $\l_1=\m_2 < 1/4$ (cf. \ref{MSM}) the
set of semistable morphisms is $W^{ss}(\sigma)=W_1$. The open subset of stable
points $W^{s}(\sigma) \subset W^{ss}(\sigma)$ is given by the conditions
$\f_{12}\nmid \f_{11}$ and $\f_{12} \nmid \f_{22}$. The geometric quotient
$W^s(\sigma)/G$ is an open subset of $W_1/\!\!/G$. 

Let $X_1$ be the locally closed subset of $\M(4,4)$ given by the relation 
$h^0(\F(-1))=1$. We equip $X_1$ with the canonical induced reduced structure.
The morphism $\rho: W_1 \to X_1$, $\rho(\f)= [\Coker(\f)]$, is surjective and
$G$-invariant. Since $W_1$ is irreducible, $X_1$ is irreducible, too. From
\ref{5.C} we know that $\rho(\f)$ is the isomorphism class of a stable sheaf if
and only if $\f$ is in $W^s(\sigma)$. Arguments similar to those in the
beginning of \ref{open2} show that $\rho(\f_1)=\rho(\f_2)$ if and only if
$\overline{G\f_1} \cap \overline{G\f_2} \neq \emptyset$. If $\F$ is stable, then
$\rho^{-1}([\F])$ is a $G$-orbit, so it has dimension equal to $\text{dim}(G)$.
Indeed, the stabilizer for any stable morphism consists only of the neutral
element of $G$. We have
\bdm
\text{dim}(X_1) = \text{dim}(W^s(\sigma)) - \text{dim}(G) = 25-9 =16.
\edm

\SEP

\begin{subsub}\label{5.D}{\bf Theorem : }
The good quotient $W_1/\!\!/G$ is isomorphic to $X_1$.
In particular, $X_1$ is a closed hypersurface of $\MM(4,4)$. The closed 
subvariety $(W_1 \setminus W^s(\sigma))/\!\!/G$ of $X_1$ is isomorphic to the
subvariety of $\MM(4,4)$ given by non stable sheaves which have a
factor of the form $\O_C(1)$ in their Jordan-H\"older filtration, for a cubic
curve $C \subset \P^2$, and is isomorphic to $\P(S^3V^*) \times \P(V^*)$.
\end{subsub}

\Proof As in the proof of \ref{3.B}, we need to construct resolution \ref{5.C}
starting from the Beilinson spectral sequence for $\F$. Tableau (2.2.3) takes
the form
\bdm
\xymatrix
{
\O(-2) & 0 & 0 \\
\O(-2) \ar[r]^{\f_3} & 4\O(-1) \ar[r]^{\f_4} & 4\O
}.
\edm
The exact sequence (2.2.5) becomes
\bdm
0 \lra \O(-2) \stackrel{\f_5}{\lra} \Coker(\f_4) \lra \F \lra 0.
\edm
Note that $\f_5$ lifts to a morphism $\O(-2) \lra 4\O$ because
$\Ext^1(\O(-2),\Coker(\f_3))=\lbrace 0\rbrace$. Combining with the exact
sequence (2.2.4) we obtain the resolution
\[\xymatrix{
0\ar[r] & \O(-2)\ar[r]^-{\begin{pmatrix}0 \\ \psi\end{pmatrix}} &
\O(-2) \oplus 4\O(-1)\ar[r]^-f & 4\O\ar[r] & \F\ar[r] & 0.
}\]
The argument from remark 6.3 in \cite{maican} shows that, up to equivalence,
\bdm
\psi^{\text{T}} = 
\begin{pmatrix}
0 & X & Y & Z
\end{pmatrix}
\edm
Now we distinguish two possibilities: firstly, up to equivalence,
\bdm
\f =
\begin{pmatrix}
\f_{11} & 0 \\
\f_{21} & \f_{22}
\end{pmatrix}
\quad \text{for a morphism} \quad \f_{11}:\O(-2) \oplus \O(-1) \lra 2\O.
\edm
We see that $\F$ surjects onto $\Coker(\f_{11})$, hence this sheaf is supported
on a proper closed subset of $\P^2$, hence $\f_{11}$ is injective, hence
$\Coker(\f_{11})$ has Hilbert polynomial $P(t)=3t+2$. This contradicts the
semistability of $\F$. The second possibility, in fact the only feasible one, is
that, modulo equivalence,
\bdm
\f = 
\begin{pmatrix}
\f_{11} & 0 \\
\f_{21} & \f_{22}
\end{pmatrix}
\quad \text{with} \quad \f_{11} =
\begin{pmatrix}
q & \ell
\end{pmatrix}
\quad \text{and} \quad \f_{22} =
\begin{pmatrix}
-Y & X & 0 \ \\
-Z & 0 \ & X \\
0 \ & -Y & Z
\end{pmatrix}
,
\edm
\bdm
\f_{11}: \O(-2) \oplus \O(-1) \lra \O, \qquad \f_{22}: 3\O(-1) \lra 3\O.
\edm
From the semistability of $\F$ we see that $\ell \neq 0$. It is also easy to see
that
$\F$ is properly semistable if and only if $\ell$ divides $q$. Let $E$ denote 
the set parametrized by $(\psi,\f)$. We can identify $E$ with an open subset
inside the affine space parametrized by the entries of $\f_{11}$ and $\f_{21}$,
so $E$ is irreducible and smooth. The subset given by the condition that $\ell$
divide $q$ has codimension 3, hence the morphism \ $\xi : E \to W_1$ \ can be
defined on its complement and then extended algebraically to $E$. Thus, for the
purpose of constructing $\xi$, we may assume in the sequel that $\ell$ does not
divide $q$. The snake lemma gives the exact sequence
\bdm
0 \lra \Ker(\f_{22}) \lra \O(-2) \lra \Ker(\f_{11}) \lra \Coker(\f_{22}) \lra \F
\lra
\edm
$$\quad\quad\quad\quad\quad\lra \Coker(\f_{11})\lra 0 .$$
As $\Ker(\f_{22})= \O(-2)$, $\Ker(\f_{11})=\O(-3)$, $\Coker(\f_{22})= \O(1)$, we
obtain
the extension
\bdm
0 \lra \O_C(1) \lra \F \lra \Coker(\f_{11}) \lra 0.
\edm
Here $C$ is the quartic curve given as the zero-set of the polynomial
\bdm
f =
\begin{pmatrix}
Z & -Y & X
\end{pmatrix}
\ \f_{21} \
\begin{pmatrix}
-\ell \\ q
\end{pmatrix}
.
\edm
From here on we construct a morphism to $W_1$ in the same manner as in the proof
of \ref{5.C}. \carre

\SEP

\begin{subsub}\label{5.E}{\bf Remark: }\rm
In the course of the above proof we have rediscovered lemma 4.10
from \cite{lepotier} which states that every stable sheaf $\F$ satisfying 
$h^0(\F(-1))=1$ occurs as an extension
\bdm
0 \lra \O_C(1) \lra \F \lra \O_S \lra 0
\edm
with $C$ a quartic curve and $S$ a zero-dimensional scheme of length 2.
The ideal of $S$ is generated by $q$ and $\ell$.
\end{subsub}

\SEP

\begin{subsub}\label{5.F}{\bf Theorem : }
There are no sheaves $\F$ in $\MM(4,4)$ satisfying the
condition \hfil\break$h^0(\F(-1)) \ge 2$.
\end{subsub}

\Proof
Let $r$ be an integer such that $1\leq r\leq 3$. It is easy, using the
descriptions of $\M(r,r)$ given in the Introduction, to see that there are no
sheaves $\E$ in $\M(r,r)$ such that \ $h^0(\E(-1))\geq 2$. It follows that if
$\F$ is semistable non stable sheaf with Hilbert polynomial $4m+4$ then we have
\ $h^0(\F(-1))\leq 1$.

According to \ref{open4}, $\M(4,4)\backslash X_0$ is the exceptional divisor of
the blowing-up \hfil\break$\M(4,4)\to N(3,4,4)$ \ along $\P^2$. Hence it is an
irreducible hypersurface. Since this hypersurface contains $X_1$ we have \
$\M(4,4)\backslash X_0=X_1$. The result follows immediately.
\carre

\end{sub}

\end{document}